\documentclass[reqno]{amsart}

\usepackage{amsmath,amssymb,amsthm,amsfonts,
	enumerate,hyperref,cleveref,bm,mleftright}

\usepackage{mathtools,leftindex,tensor,mhchem,cancel,xcolor,etoolbox,lipsum}

\newcommand{\CancelTo}[3][]{%
	\ifblank{#1}{}{%
		\renewcommand{\CancelColor}{#1}%
	}
	\cancelto{#2}{#3}%
}

\hypersetup{
	colorlinks=true,
	linkcolor=blue,
	filecolor=magenta,      
	urlcolor=cyan,
	pdftitle={Overleaf Example},
	pdfpagemode=FullScreen,
}
\usepackage{dsfont}
\usepackage[all]{xy}
\usepackage[T1]{fontenc}
\usepackage{tikz}
\usepackage{leftidx}
\usepackage{pgfplots}
\pgfplotsset{compat=1.15}
\usepackage{mathrsfs}
\usetikzlibrary{arrows}
\usepackage{isotope}

\theoremstyle{plain}
\newtheorem{thrm}{Theorem}[section]
\newtheorem{cor}[thrm]{Corollary}
\newtheorem{prop}[thrm]{Proposition}
\newtheorem{lem}[thrm]{Lemma}

\theoremstyle{definition}

\newtheorem{rem}[thrm]{Remark}
\newtheorem{exm}[thrm]{Example}

\crefrangeformat{equation}{#3(#1)#4--#5(#2)#6}

\crefname{thrm}{Theorem}{Theorems}
\crefname{theorem}{Theorem}{Theorems}
\crefname{lem}{Lemma}{Lemmas}
\crefname{cor}{Corollary}{Corollaries}
\crefname{prop}{Proposition}{Propositions}
\crefname{defn}{Definition}{Definitions}
\crefname{exm}{Example}{Examples}
\crefname{rem}{Remark}{Remarks}
\crefname{conj}{Conjecture}{Conjectures}
\crefname{quest}{Question}{Questions}
\crefname{section}{Section}{Sections}
\crefname{equation}{\unskip}{\unskip}
\crefname{enumi}{\unskip}{\unskip}
\crefname{subsection}{Subsection}{Subsections}

\newcommand{\CC}{\mathbb{C}}

\def\notop{\mathrel{\not\mskip-\thinmuskip\top}}

\begin{document}
	
\title[Additive preservers of truncations]{Characterization of minimal tripotents via annihilators and its application to the study of additive preservers of truncations}

\author[L. Li]{Lei Li}
\address[L. Li]{School of Mathematical Sciences and LPMC, Nankai University, 300071 Tianjin, China.}
\email{leilee@nankai.edu.cn}

\author[S. Liu]{Siyu Liu}
\address[S. Liu]{School of Mathematical Sciences and LPMC, Nankai University, 300071 Tianjin, China.}
\email{760659676@qq.com}

\author[A.M. Peralta]{Antonio M. Peralta}
\address[A.M. Peralta]{Instituto de Matem{\'a}ticas de la Universidad de Granada (IMAG), Departamento de An{\'a}lisis Matem{\'a}tico, Facultad de
	Ciencias, Universidad de Granada, 18071 Granada, Spain.}
\email{aperalta@ugr.es}

	
\begin{abstract} The contributions in this note begin with a new characterization of (positive) scalar multiples of minimal tripotents in a general JB$^*$-triple $E$, proving that a non-zero element $a\in E$ is a positive scalar multiple of a minimal tripotent in $E$ if, and only if, its inner quadratic annihilator (that is, the set $\leftindex^{\perp_{q}}\{a\} = \{ b\in E: \{a,b,a\} =0\}$) is maximal among all inner quadratic annihilators of single elements in $E$. We subsequently apply this characterization to the study of surjective additive maps between atomic JBW$^*$-triples preserving truncations in both directions. Let $A: E\to F$ be a surjective additive mapping between atomic JBW$^*$-triples, where $E$ contains no one-dimensional Cartan factors as direct summands. We show that $A$ preserves truncations in both directions if, and only if, there exists a  bijection $\sigma: \Gamma_1\to \Gamma_2$, a bounded family $(\gamma_k)_{k\in \Gamma_1}\subseteq \mathbb{R}^+$, and a family $(\Phi_k)_{k\in \Gamma_1},$ where each $\Phi_k$ is a (complex) linear or a conjugate-linear (isometric) triple isomorphism from $C_k$ onto $\widetilde{C}_{\sigma(k)}$ satisfying $\inf_{k} \{\gamma_k \} >0,$ and $$A(x) = \Big( \gamma_{k} \Phi_k \left(\pi_k(x)\right) \Big)_{k\in\Gamma_1},\  \hbox{ for all } x\in E,$$ where $\pi_k$ denotes the canonical projection of $E$ onto $C_k.$   
\end{abstract}
	
	\maketitle
	
	
\section{Introduction}

Preservers of truncations are currently studied in different settings, from preservers on the space $B(H)$, of all bounded linear operators on a Hilbert space $H$ (cf. \cite{Jia_Shi_Ji_AnnFunctAnn_2022,MaoJi2024,Yao_Ji_JMathResAppl_2022}), to the wider setting of Cartan factors and general JB$^*$-triples (see \cite{GarLiPeSu}). We recall that, for $a,b\in B(H)$, we say that $a$ is a truncation of $b$ if $aa^* a = ab^* a$. It is widely known that every C$^*$-algebra, $A,$ is a JB$^*$-triple with respect to the triple product \begin{equation}\label{eq Cstar triple product} \{a,b,c\} := \frac12 (a b^* c + c b^* a) \ \ (a,b,c\in A).
\end{equation} Furthermore, in terms of the just defined triple product, $a$ is a truncation of $b$ if and only if $\{a,a,a\} =\{a,b,a\}$. So, the relation ``being a truncation of'' fully makes sense in the wider setting of JB$^*$-triples (see \Cref{sec quadratic annihilators and minimal tripotents} for the detailed definitions).\smallskip

A very recent study shows that if $\Delta :A \to B$ is a {\rm(}non-necessarily linear nor continuous{\rm)} bijection between atomic JBW$^*$-triples preserving the truncation of triple products in both directions, that is, $$\begin{aligned}
	\boxed{a \mbox{ is a truncation of } \{b,c,b\}} \Leftrightarrow \boxed{\Delta(a)  \mbox{ is a truncation of } \{\Delta(b),\Delta(c),\Delta(b)\}},
\end{aligned}$$ and satisfying its restriction to each rank-one Cartan factor in $A$, if any, is a continuous mapping, then $\Delta$ is an isometric real-linear triple isomorphism \cite{GarLiPeSu}. A related precedent in the particular case that $A = B = B(H)$ was considered by X. Jia, W. Shi, and G. Ji in \cite{Jia_Shi_Ji_AnnFunctAnn_2022}. \smallskip

Just like the studies on (non-necessarily linear) bijections preserving $\lambda$-Aluthge transforms on products in \cite{Chabb2017,ChabbMbekhta2017,EssPe2018} was preceded by the description of the bijective linear transformations between von Neumann factors that commute with the $\lambda$-Aluthge transform in \cite{BoMolNag2016}, it seems natural to pose the challenge of describing all linear surjection between (atomic) JB$^*$-triples preserving truncations in both directions. Observe that this problem is independent from the commented result on bijections preserving truncations of triple products in both directions.  Moreover, it could be also asked whether the linearity of the mapping is absolutely necessary to obtain the description. Concerning this question, in \cite[Theorem 2.3]{Yao_Ji_JMathResAppl_2022}, J.~Yao and G.~Ji. prove that for each complex Hilbert space $H$ with dim$(H)\geq2$, the following statements are equivalent for every additive and surjective map $A : B(H) \to  B(H)$:\begin{enumerate}[$(a)$] 
	\item $A$ preserves truncations of operators in both directions;
	\item There exist a nonzero scalar $\alpha \in \mathbb{C}$ and maps $u,v:H\to H$ which are both unitary operators or both anti-unitary operators such that $A(x) = \alpha u x v$ for all $x\in B(H)$ or $A(x) = \alpha u x^* v$ for all $x\in B(H)$.	
\end{enumerate} Observe that statement $(b)$ above is equivalent to say that $A$ is a positive multiple of a (complex) linear or a conjugate-linear triple automorphism on $B(H)$.\smallskip 

As commented before $B(H)$ is an example of C$^*$-algebra, and a JB$^*$-triple. It can be also regarded inside the collection of complex Banach spaces called Cartan factors. We recall their definition with the aim to provide some basic background on Cartan factors for the readers. Let $H$ and $K$ be complex Hilbert spaces. The spaces $B(H, K)$ of all bounded linear operators between $H$ and $K$ form the so-called \textit{Cartan factors of type $1$}. Suppose next that $j: H\to H$ is a conjugation on $H$ (i.e., a conjugate-linear isometry of period $2$). The subspaces  of $B(H)$ given by $C_2= \{a\in B(H): a=-ja^*j\}$ and $C_3=\{a\in B(H): a=ja^*j\}$ determine the classes of \textit{Cartan factors of type $2$} and $3$, respectively. All these spaces are equipped with the triple product given in \eqref{eq Cstar triple product}. A Banach space $V$ is called a \textit{Cartan factor of type $4$} or a \textit{spin factor} if it admits a complete inner product $\langle \cdot|\cdot\rangle$ and a conjugation $x\mapsto \bar{x}$, for which the norm of $V$ is given by
\[\|x\|^2=\langle x|x \rangle +\sqrt{\langle x|x\rangle^2-|\langle x|\bar{x}\rangle|^2}\] and the triple product of $V$ is defined by $$\{x, y, z\} = \langle x|y\rangle z + \langle z|y\rangle  x -\langle x|\overline{z}\rangle \overline{y}, \ \ (x,y,z\in V).$$ The \textit{Cartan factor of type $5$ and $6$} (also called  \textit{exceptional Cartan factors}) are the finite dimensional Banach spaces $M_{1,2}(\mathbb{O})$ and $H_3(\mathbb{O}),$ of all $1\times 2$ matrices with entries in the (complex) octonions $\mathbb{O}$, and all $3\times 3$ hermitian matrices with entries in $\mathbb{O}$, respectively (see  \cite{HamKalPe2023} for more details). \smallskip
     
In our main result in this note we study surjective additive mappings preserving truncations in both directions between JB$^*$-triples which can be written as $\ell_{\infty}$-sums of families of Cartan factors. Our main conclusion, established in \Cref{t main theorem}, reads as follows: Let $\displaystyle E = \bigoplus_{k\in \Gamma_1}^{\ell_{\infty}} C_k$ and $\displaystyle F = \bigoplus_{{j\in \Gamma_2}}^{\ell_{\infty}} \widetilde{C}_j$  be atomic JBW$^*$-triples,  where $\{C_k\}_{k\in \Gamma_1}$ and $\{\widetilde{C}_j\}_{j\in \Gamma_2}$ are two families of Cartan factors. Suppose, additionally, that $E$ contains no one-dimensional Cartan factors as direct summands (i.e., dim$(C_k)\geq 2$ for all $k$). Let $A: E\to F$ be a surjective additive mapping. Then the following statements are equivalent:  
\begin{enumerate}[$(a)$]
	\item $A$ preserves truncations in both directions. 
	\item There exists a  bijection $\sigma: \Gamma_1\to \Gamma_2$, a bounded family $(\gamma_k)_{k\in \Gamma_1}\subseteq \mathbb{R}^+$, and a family $(\Phi_k)_{k\in \Gamma_1},$ where each $\Phi_k$ is a (complex) linear or a conjugate-linear (isometric) triple isomorphism from $C_k$ onto $\widetilde{C}_{\sigma(k)}$ satisfying $\inf_{k} \{\gamma_k \} >0,$ and $$A(x) = \Big( \gamma_{k} \Phi_k \left(\pi_k(x)\right) \Big)_{k\in\Gamma_1},\  \hbox{ for all } x\in E,$$ where $\pi_k$ denotes the canonical projection of $E$ onto $C_k.$  
\end{enumerate} In particular, every surjective additive mapping $A: E\to F$ preserving truncations in both directions is a real-linear bijection, and the JB$^*$-triples $E$ and $F$ are isometrically triple isomorphic. Moreover, if $A$ satisfies $(a)$ or $(b)$ and there exists an element $x_0$ in $E$ satisfying $\|A\pi_k (x_0)\| = \|\pi_k (x_0)\| \neq 0$, for all $k\in \Gamma_1$, the mapping $A$ is a real-linear {\rm(}isometric{\rm)} triple isomorphism. Note that the implication $(b)\Rightarrow (a)$ is clear. We shall see in \Cref{counterexample for C complex numbers} that the hypothesis affirming that dim$(C_k)\geq 2$ for all $k$ cannot be relaxed.\smallskip

It is known that a JB$^*$-triple has rank-one if and only if it is a complex Hilbert space regarded as a type $1$ Cartan factor (cf. \Cref{sec: rank-one JB*-triples}). Therefore, every additive mapping $A$ between two rank-one JB$^*$-triples is nothing but an additive mapping between two Hilbert spaces. We shall see below that, in such a case, $A$ preserves truncations if, and only if, it preserves Euclidean orthogonality between elements in the corresponding Hilbert spaces (cf. \Cref{r quadratic annihilator and Euclidean orthogonality}). We are thus in a position to apply a recent result on additive preservers of Euclidean orthogonality between complex Hilbert spaces, obtained in \cite{LiLiuPe24}, to deduce that if $A$ is a surjective additive mapping between rank-one JB$^*$-triples of dimension $\geq 2,$ and $A$ preserves truncations (in one direction), then $A$ is a positive scalar multiple of a {\rm(}complex{\rm)} linear or a conjugate-linear isometry (see \Cref{p trunc preservers between Hilbert spaces}).\smallskip

The proof of the main result in this note has also required the development of a tool of independent interest. In \Cref{th characterization of minimal tripotents} we establish that a non-zero element $a$ in a JB$^*$-triple $E$ is a {\rm(}positive{\rm)} scalar multiple of a minimal tripotent in $E$, if and only if, its inner quadratic annihilator $\leftindex^{\perp_{q}}\{a\} = \{ b\in E: \{a,b,a\} =0\}$ is maximal among all inner quadratic annihilators of single elements in $E$. This result generalises a previous characterization of rank-one operators in $B(H)$ obtained by J. Yao and G. Ji in \cite[Lemma 2.2]{Yao_Ji_JMathResAppl_2022}. \smallskip

We have essentially described the main results in sections 2 and 3. \Cref{sec: general properties of truncation preservers on general JB*-triples} is devoted to analyse the basic properties of the surjective additive maps preserving truncations in both directions between two general JB$^*$-triples $E$ and $F$. It is shown that any such a mapping must be automatically injective and maps scalar multiples of minimal tripotents in $E$ (if any) to scalar multiples of minimal tripotents in $F$. Moreover, for each minimal tripotent $e$ in $E$ we have $A (\mathbb{C} e) = \mathbb{C} A(e) = \mathbb{C} r(A(e))$, where $r(A(e))$ denotes the range tripotent of $A(e)$, which in this case lies in $F$ (see \Cref{l basic properties of additive preservers of truncations into abd images of min trip}). Further technical conclusions assure that if $w_1$ and $w_2$ are two collinear (respectively, orthogonal) minimal tripotents in $E$, the corresponding range tripotents $r(A(w_1))$ and $r(A(w_2))$ are minimal and collinear (respectively, orthogonal) in $F$ (cf. Propositions~\ref{prop preservation of collinearity} and \ref{prop preservation of orthogonality}).\smallskip

The concluding section is devoted to establish our main result. We prove that if $A: E\to F$ is a surjective additive mapping preserving truncations in both directions between two atomic JBW$^*$-triples, where $E$ contains no one-dimensional Cartan factors as direct summands, the mapping $A_r : e\mapsto r(A(r))$, sending each minimal tripotent in $E$ to the range tripotent of $A(e)$ in $F$ is a bijective map preserving collinearity and orthogonality in both directions from the set $\mathcal{U}_{min} (E)$ of all minimal tripotents in $E$ onto the corresponding set $\mathcal{U}_{min} (F)$ (see \Cref{c each Cartan factor is mapped to a Cartan factor 2 Badalin-3}). Furthermore, when restricted to each Cartan factor $C_k$ in the decomposition of $E$, the mapping $A_r$ either preserves or reverses triple transition pseudo-probabilities between all minimal tripotents in $C_k$, actually $A|_{C_k}$ must be a positive scalar multiple of a linear or conjugate-linear triple isomorphism from $C_k$ onto a Cartan factor in the decomposition of $F$ (see \Cref{c each Cartan factor is mapped to a Cartan factor 2 Badalin-5}). We can finally apply the description of all bijective maps preserving triple transition pseudo-probabilities between the sets of minimal tripotents in two atomic JBW$^*$-triples borrowed from \cite{Peralta_ResMath_2023} to complete the proof of our main result.

\section{Quadratic annihilators and minimal tripotents}\label{sec quadratic annihilators and minimal tripotents}

It is perhaps worth to begin this section by recalling the definition of JB$^*$-triples.  A \emph{JB$^*$-triple} (see \cite{Kaup_RiemanMap}) is a Banach space together with a continuous triple product $$ \{\cdot, \cdot, \cdot  \}: E \times E \times E \to E : (x, y,z ) \mapsto \{ x,y,z \}$$
for all $x,y,z \in E$, which is linear in the first and third variables, conjugate-linear in the middle one, and satisfies the following axioms: 
\begin{enumerate}[$(i)$]
	\item\label{Jordan identity} For each $a,b$ in $E$, the operators $L(a,b): E\to E$, $L(a,b)(c) := \{a,b,c\}$ satisfy the identity
	$$L(w,v)\{x,y,z \} = \{L(w,v)x,y,z \} - \{x,L(v,w)y,z \}+\{x,y, L(w,v)z\},$$ for all $x,y,z,w,v \in E$; \hfill ($\emph{Jordan identity}$)  
	\item For each $x\in E$, the operator $L(x,x)$ is hermitian with non-negative spectrum;
	\item $\| \{x,x,x\}  \| = \| x\|^3$ for all $x \in E$.\hfill ($\emph{Gelfand-Naimark axiom}$) 
\end{enumerate} All Cartan factors are JB$^*$-triples.\smallskip 

A JBW$^*$-triple is a JB$^*$-triple whose underlying Banach space is a dual Banach space. Every JBW$^*$-triple admits an unique (isometric) predual and its triple product is separately continuous \cite{Bar_Tim_MathScand_1986}. It is worth to note that all von Neumann algebras are examples of JBW$^*$-triples, and the bidual, $E^{**},$ of a JB$^*$-triple $E$ is a JBW$^*$-triple \cite{Dineen_theseconddual_1986}. JBW$^*$-triples satisfy extra geometric algebraic properties, for example, they contain an abundant collection of tripotents; since tripotents in a JB$^*$-triple $E$ coincide with the extreme points of the closed unit ball of $E$ (see, for example, \cite[Corollary 4.8]{EDW_RUTT_JLMS_1988} or \cite[Theorem 4.2.34]{CabreraPalaciosBook}). Note that the set of all tripotents in a general JB$^*$-triple might be empty.\smallskip

The relation ``being a truncation of'' between elements in a JB$^*$-triple can be better understood via inner quadratic annihilators. The \emph{inner quadratic annihilator} of a subset $S$ of a JB$^*$-triple $E$ is defined as the set $$\leftindex^{\perp_{q}}S=\{ b\in E: Q(s)(b) =0, \ \forall\, s\in S\} = \bigcap_{s\in S} \ker(Q(s)).$$ We note that, given an element $s$ in a JB$^*$-triple $E$, we denote by $Q(s)$ the conjugate-linear mapping on $E$ given by $Q(s) (x) = \{s,x,s\}.$ Observe that $\leftindex^{\perp_{q}}S$ is a closed linear subspace of $E$ since $Q(s)$ is a bounded conjugate-linear operator on $E$ for every $s$. It is known that for each element $a$ in $E$ we have \begin{equation}\label{eq chracterization of inner quadr annihilator of an element} \leftindex^{\perp_{q}}\{a\} = E\cap \Big( E_0^{**}(r(a)) \oplus E_1^{**}(r(a))\Big) = \{ x\in E: P_2(r(a)) (x) =0\},
\end{equation} where $r(a)$ denotes the range tripotent of $a$ in $E^{**}$, and $E_j(e)$ stands for the Peirce-$j$ subspace associated to $e$, and $P_j (e)$ stands for the projection  of $E$ onto $E_j (e)$ (cf. \cite[Lemma 2.2]{GarLiPeSu}). In case that $E$ is a JBW$^*$-triple, we can further conclude that $$ \leftindex^{\perp_{q}}\{a\} =  E_0(r(a)) \oplus E_1(r(a)) = \{ x\in E: P_2(r(a)) (x) =0\},$$ where $r(a)$ denotes the range tripotent of $a$ in $E.$ \smallskip

Some of the notions employed in the previous paragraph are perhaps less known for non-experts. The missing details will be fully explained now. For example, an element $e$ in a JB$^*$-triple $E$ is called a \emph{tripotent} if $\{e,e,e\} =e$. In such a case, we can always $E$ according to the so-called \emph{Peirce decomposition} induced by $e$ in the form
\begin{equation}
E = E_{0}(e) \oplus E_{1}(e) \oplus E_{2}(e), 
\end{equation}
where each $E_{j}(e) := \{  a\in E: L(e,e)(a) = \frac{j}{2} a \}$ is a subtriple of $E$ called the \emph{Peirce $j$-subspace} ($j = 0,1,2$). Triple products among elements in Peirce subspaces follow certain patterns known as \emph{Peirce rules} or  \emph{Peirce arithmetics}, namely, $\{ E_{i}(e),E_{j}(e),E_{k}(e) \} \subseteq E_{i-j+k}(e)$ for all $i,j,k = \{0,1,2\}$ and 
$$ \{E_{0}(e), E_{2}(e), E \} = \{ E_{2}(e), E_{0}(e),E \}  = 0,$$
where $E_{i-j+k}(e) = \{0\} $ if $i-j+k \neq  \{0,1,2\}$. The Peirce 2-subspace $E_{2}(e)$ is actually a unital JB$^*$-algebra with identity $e$, where the Jordan product and involution operations are given by $a \circ_{e} b := \{a,e,b \}$ and $a ^{*_{e}} := \{e,a,e \}$, respectively (cf. \cite{Horn_MathScand_1987} or \cite[Fact 4.2.14, Proposition 4.2.22, and Corollary 4.2.30]{CabreraPalaciosBook}). An illustrative example is given by C$^*$-algebras (and hence by $B(H)$ spaces), regarded as JB$^*$-triples with respect to the triple product in \eqref{eq Cstar triple product}, where tripotents correspond to partial isometries.\smallskip

We can classify tripotents in terms of the summands in the corresponding Peirce decomposition. For example, a tripotent $e$ is called (\emph{algebraically}) \emph{minimal} (respectively, \emph{complete} or \emph{maximal}) if $E_{2}(e) = \CC e \neq \{ 0\}$ (respectively, $E_{0}(e) = \{ 0\}$). We shall write $\mathcal{U}(E)$,  $\mathcal{U}_{min}(E)$ and $\mathcal{U}_{max}(E)$ for the sets of all tripotents, all minimal tripotents, and maximal tripotents in $E$, respectively. We employed the word ``algebraic'' because there is also a notion of minimality associated with the partial ordering on tripotents.\smallskip

The characterization of the inner quadratic annihilator in \Cref{eq chracterization of inner quadr annihilator of an element} can be now applied to deduce the following result.

\begin{lem}\label{l lineraly independence of multiples of minimal tripotents} Let $e,v$ be two minimal tripotents in a JB$^*$-triple $E$, and let us pick two non-zero elements $a\in \mathbb{C} e$ and $b\in \mathbb{C} v$. Then $a$ and $b$ are linearly dependent if, and only if, $ \leftindex^{\perp_{q}}\{a\} =  \leftindex^{\perp_{q}}\{b\}.$  
\end{lem}

Elements $x,y$ in a JB$^*$-triple $E$ are said to be \emph{orthogonal} ($x \perp y$ in short) if $L(x,y) = 0$ (see, for example, Lemma 1 and the comments around it in \cite{BFPGMP08} for additional properties). When the relation of orthogonality is restricted to $\mathcal{U}(E)$, it is known that $w,v\in \mathcal{U}(E)$, are orthogonal if, and only if, $w \in E_{0}(v)$. It is known that $a\perp b$ in $E$ implies that they are $M$-orthogonal, that is, $\|a+b\| = \max\{\|a\|, \|b\|\}$ (see \cite[Lemma 1.3$(a)$]{FriedRusso_Predual}).\smallskip 

Let us recall that two tripotents $w_1,w_2$ in a JB$^*$-triple $E$ are called \emph{collinear} ($w_1\top w_2$ in short) if $w_i \in E_1 (w_j)$ for all $i\neq j$ in $\{1,2\}$. We say that $w_2$ \emph{governs} $w_1$ ($w_2 \vdash w_1$ in short) if $w_1\in E_{2} (w_2)$ and $w_2\in E_{1} (w_1)$. \smallskip

The natural partial ordering on $\mathcal{U}(E)$ is defined as follows: for $e,u\in \mathcal{U}(E)$ the symbol $e \leq u$ means that $u-e \in \mathcal{U}(E) $ and $u -e \perp e$, or equivalently, $e$ is a projection in the unital JB$^*$-algebra $E_{2}(u)$. We note that these relations of orthogonality and partial ordering among elements in $\mathcal{U}(E)$ agree with the original notions of orthogonality and order for C$^*$-algebras in that particular setting. Clearly, every (algebraic) minimal tripotent is order minimal, however the reciprocal statement is not necessarily true as shown by the unit element in $C[0,1]$. Under the stronger assumption that $E$ is a JBW$^*$-triple minimal and order minimal tripotents agree (cf. \cite[Corollary 4.8]{EDW_RUTT_JLMS_1988} and \cite[Lemma 4.7]{Battaglia_1991}). \smallskip

For each element $a$ in a JBW$^*$-triple $E$, there exists a smallest tripotent $e$ in $E$ (called the \emph{range tripotent} of $a$) satisfying that $a$ is a positive element in the JBW$^*$-algebra $E_{2} (e)$. We write $r(a)$ for the range tripotent of $a$.  It is known that $r(a)$ also coincides the range projection of $a$ in $E_2 (r(a))$.

\begin{rem}\label{r chracterization of truncations} It is known from \cite[Lemma 2.4]{GarLiPeSu} that an element $a$ in a JB$^*$-triple $E$ is a truncation of another element $b$ in $E$ if, and only if, one of the following equivalent statements holds:\label{characterization of triple truncations}
	\begin{enumerate}[{$(a)$}]
		\item $b=a+z$ for some $z\in E$ with $\{a,z,a\}=0$ {\rm(}i.e. $z\in \leftindex^{\perp_{q}}\{a\} ${\rm)};
		\item $a=P_2(r(a)) (b)$ where $r(a)$ is the range tripotent of $a$ in $E^{**}$. 
	\end{enumerate}
\end{rem}

The following lemma is an obvious consequence of the characterization given in \Cref{r chracterization of truncations}.

\begin{lem}\label{l additive preservers of truncations coincide with additive preservers of inner quadratic annihilators} Let $A: E\to F$ be an additive mapping between JB$^*$-triples. Then $A$ preserves truncations if, and only if, it preserves inner quadratic annihilators. 
\end{lem}

Preservers of rank-one operators on subalgebras of the space $B(X),$ of all bounded linear operators on a nontrivial real or complex Banach space $X,$ have been intensively studied, and we have  (see, for example, \cite{OmlaSemrl93} and the historical review therein).\smallskip

We recall that a JBW$^*$-triple is called \emph{atomic} if it coincides with the weak$^*$-closure of the linear span of all its minimal tripotents (cf.\cite{FriedRusso_Predual,FriedRusso_GelfNaim}). Atomic JB$^*$-triple can be concretely represented as $\ell_{\infty}$-sums of families of Cartan factors \cite[Proposition 2 and Theorem E]{FriedRusso_GelfNaim}.\smallskip

In $B(X)$ the notion of rank-one operator is quite natural and needs no extra explanation. But, do we have an equivalent notion in more general structures? In \cite[Lemma 2.2]{Yao_Ji_JMathResAppl_2022} J. Yao and G. Ji established that a non-zero operator $a$ in $B(H)$ is a rank-one operator if, and only if, $\leftindex^{\perp_{q}}\{a\} = \{x\in B(H) : a x^* a =0 \}$ is maximal among all inner quadratic annihilators of elements in $B(H)$. It is natural to ask whether a similar conclusion holds when $B(H)$ is replaced with a von Neumann algebra, a C$^*$-algebra or a general JB$^*$-triple. Our first main result presents an argument to solve this question.

\begin{thrm}\label{th characterization of minimal tripotents} Let $a$ be a non-zero element in a JB$^*$-triple $E$. Then the following statements are equivalent.\begin{enumerate}[$(a)$]
		\item $a$ is a {\rm(}positive{\rm)} scalar multiple of a minimal tripotent in $E$.
		\item The inner quadratic annihilator of $a$, $\leftindex^{\perp_{q}}\{a\}$, is a maximal element in the set of all inner quadratic annihilators of single elements in $E$.  
	\end{enumerate}
\end{thrm}

\begin{proof}
$(a)\Rightarrow (b)$ Suppose $a= \lambda e$, where $\lambda\in \mathbb{C}\backslash \{0\}$ and $e$ is a minimal tripotent. 	It is easy to see that $\leftindex^{\perp_{q}}\{a\} = \leftindex^{\perp_{q}}\{e\} = E_1 (e) \oplus E_0 (e)$ (cf. \eqref{eq chracterization of inner quadr annihilator of an element}, for latter equality). Since $E= E_2(e) \oplus E_1 (e) \oplus E_0 (e)$ with $E_2 (e) = \mathbb{C} e$, it follows that $\leftindex^{\perp_{q}}\{a\}$ is a one-codimensional closed subspace of $E$, and hence maximal among subspaces. Having in mind that the inner quadratic annihilator associated to each subset of $E$ is a closed linear subspace, we deduce the desired property for $\leftindex^{\perp_{q}}\{a\}$.\smallskip
	
$(b)\Rightarrow (a)$ As we shall see now, it is relatively easy to prove that $a$ is a positive scalar multiple of a tripotent, but showing that this tripotent is in fact minimal in $E$ is more difficult. By local theory, the JB$^*$-subtriple of $E$ generated by the element $a$, denoted by $E_a$,  is (isometrically) JB$^*$-triple isomorphic to a commutative C$^*$-algebra of the form $C_0(\Omega)$ for some locally compact Hausdorff space $\Omega \subseteq (0,\|a\|],$ such that $\Omega\cup \{0\}$ is compact, and under the corresponding identification $a$ is identified with the inclusion mapping of $\Omega$ into $\mathbb{C}$ (cf.
\cite[Corollary 4.8]{Kaup77MathAnn}, \cite[Corollary 1.15]{Kaup_RiemanMap} and
\cite{FriedmanRusso82TAMS,FriRuCommutative}). It is also known that the range tripotent of $a$ in $E^{**}$ is a unit element in $E_a^{**}$. If $\Omega$ contains more than one point, we can take two orthogonal non-zero elements $c,b$ in $E_a \cong  C_0(\Omega)$  (positive in the local order) whose range tripotents satisfy $r(b), r(c)\leq r(a)$ in $E_a^{**},$ and also as tripotents in $E^{**}$. It follows from the inequality $r(b)\leq r(a)$ and \eqref{eq chracterization of inner quadr annihilator of an element} that $\leftindex^{\perp_{q}}\{a\} \subseteq \leftindex^{\perp_{q}}\{b\}$. Since $c\perp b$ we have $c\in  \leftindex^{\perp_{q}}\{b\},$ and $\{a,c,a\}\neq 0$ gives $c\notin \leftindex^{\perp_{q}}\{a\}$. Therefore $\leftindex^{\perp_{q}}\{a\}\subsetneqq \leftindex^{\perp_{q}}\{b\},$ which contradicts the maximality of $\leftindex^{\perp_{q}}\{a\}$. Therefore $\Omega$ reduces to a single element, and hence $a$ is a positive scalar multiple of a tripotent $e$ in $E_a,$ and clearly a tripotent in $E$. We have thus proved that $a = \lambda e$, for $\lambda\in \mathbb{R}^{+}$ and a tripotent $e$ in $E$.\smallskip

Our next goal will consist in proving that $e$ is a minimal tripotent in $E$. We recall that we can decompose $E^{**}$ as the direct sum of two weak$^*$-closed ideal $A$ and $N$, where $A,$ called the atomic part of $E^{**}$, is generated by all minimal tripotent in $E^{**}$, while $N$ contains no minimal tripotents (cf. \cite[Theorem 2]{FriedRusso_Predual}). It is also known that if $\pi: E^{**}\to A$ denotes the canonical projection of $E^{**}$ onto $A$, and $\iota_{E}: E\hookrightarrow E^{**}$ stands for the canonical inclusion, the mapping $\Phi = \pi\circ \iota_{E}: E\to A$ is an isometric triple monomorphism with weak$^*$-dense range (see \cite[Proposition 1 and its proof]{FriedRusso_GelfNaim}). We further know from the just quoted reference that $A$ is an atomic JBW$^*$-triple, that is, it can be written as an $\ell_{\infty}$-sum of Cartan factors (cf. \cite[Proposition 2 and Theorem E]{FriedRusso_GelfNaim}). It is clear that in order to prove that $e$ is minimal in $E$ it suffices to show that $\Phi (e)$ is minimal in $A$.\smallskip

Suppose that $\Phi (e)$ is not minimal in $A$. Since the latter is an atomic JBW$^*$-triple, we can find two orthogonal minimal tripotents $v_1,v_2$ in $A$ such that $\Phi (e) \geq v_1,v_2$ in $A$. Since $e\in E$, it is easy to check, from the weak$^*$-density of $\Phi(E)$ in $E$ and the weak$^*$-continuity of Peirce projections, that $\Phi(E)_2 (\Phi(e)) = A_2(e)\cap \Phi(E)$ is weak$^*$-dense in $A_2(e)$. Now, Kadison's transitivity theorem for JB$^*$-triples (\cite[Theorem 3.3]{BunFerMartPe06} or \cite[Corollary 1.3]{FerPe2007} and the explanation in \cite[\S 3]{FerPe2010} for completeness), applied to $A$ and $\Phi(E)$, implies the existence of two norm-one pairwise orthogonal elements $\Phi(b),\Phi(c)$ in $\Phi(E)_2 (\Phi(e))$ such that $\Phi (b) = v_1 + P_0 (v_1) (\Phi(b))$ and $\Phi (c) = v_2 + P_0 (v_2) (\Phi(c)).$ Since $\Phi (b) \in  \Phi(E)_2 (\Phi(e))$ it can be easily deduced from Peirce arithmetic that $ \leftindex^{\perp_{q}}\{\Phi(e)\} = \leftindex^{\perp_{q}}\{\Phi(a)\} \subseteq \leftindex^{\perp_{q}}\{\Phi(b)\}$. It follows from  $\Phi (c)\perp \Phi(b)$ that $\Phi(c)\in \leftindex^{\perp_{q}}\{\Phi(b)\}$, while $\Phi (c)\notin \leftindex^{\perp_{q}}\{\Phi(a)\}$ because $\Phi (c) \in \Phi(E)_2 (\Phi(e))$. This contradicts the assumed maximality of $\leftindex^{\perp_{q}}\{a\}$. Therefore $\Phi (e)$ must be minimal in $A$ (equivalently, $e$ must be minimal in $E$), as desired.  		
\end{proof}

 The next corollary is an interesting characterization of minimal tripotents. 
 
 \begin{cor}\label{c characterization of minimal tripotents} Let $a$ be a norm-one element in a JB$^*$-triple $E$. Then the following statements are equivalent.\begin{enumerate}[$(a)$]
 		\item $a$ is a minimal tripotent in $E$.
 		\item The inner quadratic annihilator of $a$, $\leftindex^{\perp_{q}}\{a\}$, is maximal among the set of all inner quadratic annihilators of elements in $E$.  
 	\end{enumerate}
 \end{cor} 
 
From \Cref{th characterization of minimal tripotents} one can derive a characterisation of positive scalar multiples of minimal partial isometries  in C$^*$-algebras, which seems to be also new.

\begin{cor}\label{c chract of minimal parial isometries}  Let $a$ be a non-zero element in a C$^*$-algebra $A$. Then the following statements are equivalent.\begin{enumerate}[$(a)$]
		\item $a$ is a {\rm(}positive{\rm)} scalar multiple of a minimal partial isometry in $A$.
		\item The inner quadratic annihilator of $a$, $\leftindex^{\perp_{q}}\{a\}$, is maximal among the set of all inner quadratic annihilators of elements in $A$.  
	\end{enumerate}
\end{cor}

\section{The case of rank-one JB$^*$-triples}\label{sec: rank-one JB*-triples}

A general JB$^*$-triple need not contain a single non-zero tripotent. However, every Cartan factor coincides with the weak$^*$-closure of the linear of its minimal tripotents, and the same occurs to any arbitrary orthogonal or direct sum of a family of Cartan factors. Actually a JBW$^*$-triple coincides with the weak$^*$-closure of the linear span of its minimal tripotents (also called atomic) if, and only if, it is a direct sum of a family of Cartan factors \cite[Proposition 2]{FriedRusso_GelfNaim}. It is perhaps worth to review the notion of rank. A tripotent $e$ in a JB$^*$-triple $E$ has rank--$n$ with $n\in \mathbb{N}$ (denoted by $r(e) =n$), if it can be written as the sum of $n$ mutually orthogonal minimal tripotents in $E$. A subset $S$ of $E$ is called \emph{orthogonal} if $0 \notin S$ and $x \perp y$ for every $x\neq y$ in $S$. The minimal cardinal number $r$ satisfying $\hbox{card}(S) \leq r$ for every orthogonal subset $S \subseteq E$ is called the \emph{rank} of $E$ (cf. \cite{Kaup_ManMath_1997}, \cite{BuChu92} and \cite{BeLoPeRo} for basic results on the rank of a Cartan factor and a JBW$^*$-triple, and its connection with the reflexivity of the underlying Banach space). 
\smallskip

Let us illustrate the previous notions with some examples.  Let $H$ be a complex Hilbert space regarded as a type 1 Cartan factor of the form $H\cong B(H, \mathbb{C})$, that is we consider its Hilbert norm and the triple product given by $$\{a,b,c\} = \frac12 (\langle a|  b\rangle c + \langle c | b\rangle a), \ \ (a,b,c\in H).$$ It is easy to see that every norm-one element in $H$ is a minimal and complete tripotent in $H$, and the latter has rank--one. Actually, a JB$^*$-triple has finite rank if, and only if, it is reflexive, and it has rank-one precisely when it is isometrically isomorphic to a complex Hilbert space regarded as a type 1 Cartan factor (see the discussion in \cite[\S 3]{BeLoPeRo}, \cite[page 210]{Kaup_ManMath_1997}, \cite[Corollary in page 308]{Dang_Friedm_MathScand_1987}). In case that $H$ and $K$ are complex Hilbert spaces with dimensions $m\leq n$, respectively, the JB$^*$-triple $B(H,K)$ has rank--$m$.\smallskip

The relationship ``being a trunction of'' admits a particular nice characterization in the case of rank-one JB$^*$-triples. 

\begin{rem}\label{r quadratic annihilator and Euclidean orthogonality} Let $H$ be a Hilbert space regarded as a type 1 Cartan factor, and let $x,y$ be non-zero elements in $H$. It is easy to see that $x$ is a truncation of $y$ (i.e. $\{x,x,x\}=\{x,y,x\}$) if, and only if, $\langle y|x\rangle = \langle x|x\rangle$, equivalently, $y = x + z$ with $\langle x|z\rangle =  0$. We can alternatively say that the inner quadratic annihilator of $x$ is precisely the (Euclidean) orthogonal complement of $x$, that is, $$\leftindex^{\perp_{q}} \{x\} = \{x\}^{\perp_2} =\{z\in H : \langle z|x\rangle = 0\}.$$ In this setting we write $x\perp_2 y $ if $\langle x| y\rangle =0$. So, the study of mappings between Hilbert spaces preserving (Euclidean) orthogonality is equivalent to the study of mappings preserving truncations.\smallskip

	Every non-zero element $x$ in $H$ can be written in the form $x =\|x\| \frac{x}{\|x\|}$, where $\frac{x}{\|x\|}$ is a minimal tripotent in $H$. It is easy to check that $x,y\in H\backslash \{0\}$ satisfy $x\perp_2 y$ if, and only if, $\frac{x}{\|x\|}$ and $\frac{y}{\|y\|}$ are collinear. 
\end{rem}	

Motivated by the previous remark, we recall a well-known definition. Let $H$ and $K$ be two real or complex Hilbert spaces. A mapping $\Delta : H\to K$ preserves (Euclidean) orthogonality if $$\forall x,y \in  H, \  x \perp_2 y\Rightarrow  \Delta (x) \perp_2 \Delta (y).$$  The mapping $\Delta$ preserves (Euclidean) orthogonality in both directions if the equivalence $x \perp_2 y\Leftrightarrow  \Delta (x) \perp_2 \Delta (y)$ holds for all $x,y\in H$.\smallskip

\begin{exm}\label{counterexample for C complex numbers}
In case that our JB$^*$-triple $E$ reduces to the complex field, the characterization of additive surjective maps preserving truncations in both directions on $E$ is simply hopeless, since every additive surjective mapping on $\mathbb{C}$ preserves truncations (equivalently, (Euclidean) orthogonality) in both directions. We can easily find an additive bijection on $\mathbb{C}$ preserving orthogonality in both directions (cf. the introduction of \cite{LiLiuPe24}). 
\end{exm}

We shall see next that the counterexample above can only occur when the Hilbert space one--dimensional.

\begin{prop}\label{p trunc preservers between Hilbert spaces} Let $E$ and $F$ be rank-one JB$^*$-triples with dim$(E)\geq 2$. Suppose that $A: E\to F$ is a surjective additive mapping preserving truncations, then $A$ is a positive scalar multiple of a {\rm(}complex{\rm)} linear or a conjugate-linear isometry.  
\end{prop}

\begin{proof} As we commented above, under our hypotheses, $E$ and $F$ must be complex Hilbert spaces regarded as type 1 Cartan factors. By \Cref{l additive preservers of truncations coincide with additive preservers of inner quadratic annihilators} and \Cref{r quadratic annihilator and Euclidean orthogonality}, $A$ preserves (Euclidean) orthogonality of the underlying Hilbert spaces. Since $A$ is surjective and additive by hypotheses, the desired conclusion is now a consequence of Theorem~2.1 in \cite{LiLiuPe24}.
\end{proof}

\begin{rem}\label{r chracterization of C} Let $E$ be a JB$^*$-triple. Then $$\boxed{E  \vspace*{5mm} = \mathbb{C}  \Leftrightarrow \leftindex^{\perp_{q}}\{a\} =\{0\}, \hbox{ for all } a\in E\backslash\{0\}}$$ The ``only if'' implication is clear. For the ``if'' implication, observe that if $a$ and $b$ are two orthogonal non-zero elements in $E$, we have $0\neq b\in \leftindex^{\perp_{q}}\{a\}$, which contradicts the assumption. So, $E$ has rank-one, and thus $E$ is a complex Hilbert space. If dim$(E)\geq 2$, we can take two non-zero vectors $a,b\in E$ with $\langle a | b\rangle =0$. In this case $0\neq b \in \leftindex^{\perp_{q}}\{a\}$, which is also a contradiction.    
\end{rem}

\section{Additive preservers of truncations}\label{sec: general properties of truncation preservers on general JB*-triples}

We have already observed in \Cref{l additive preservers of truncations coincide with additive preservers of inner quadratic annihilators} that an additive mapping between JB$^*$-triples preserves truncations if, and only if, it preserves inner quadratic annihilators. In this section we study the first properties of surjective additive preservers of truncations. It is well-known that a mapping between real or complex normed spaces is additive if, and only if, it is $\mathbb{Q}$-linear. We shall use this fact without any further mentioning. From now on, the symbol $\mathbb{T}$ will stand for the unit sphere of $\mathbb{C}$. 

\begin{lem}\label{l basic properties of additive preservers of truncations into abd images of min trip} Let $A: E\to F$ be a surjective additive mapping between JB$^*$-triples. Suppose additionally that $A$ preserves truncations of elements in both directions, that is, $$\{a,b,a\} =\{a,a,a\} \hbox{ in E } \Leftrightarrow \{A(a),A(b),A(a)\} =\{A(a),A(a),A(a)\}.$$ Then the following statements hold: \begin{enumerate}[$(a)$]
		\item $A$ is injective.
		\item dim$(E)\geq 2$ if, and only if, dim$(F)\geq 2$.
		\item $A$ maps minimal tripotents in $E$ to positive scalar multiples of minimal tripotents in $F$, Moreover, if $e$ is a minimal tripotent in $E$ and $\lambda$ is a non-zero complex number, then $A(\lambda e)$ is a positive scalar multiple of a minimal tripotent in $F$, that is, $A(\mathbb{C}\ \mathcal{U}_{min} (E)) = \mathbb{C}\ \mathcal{U}_{min} (F)$.  
		\item For each minimal tripotent $e$ in $E$ we have $A (\mathbb{C} e) = \mathbb{C} A(e) = \mathbb{C} r(A(e))$.  
		\item For each minimal tripotent $e\in E$, we have $r(A^{-1} (r(A(e)))) \in \mathbb{T} e$.
	\end{enumerate}
\end{lem}

\begin{proof}
$(a)$ By the additivity of $A$ it suffices to prove that $A(a)=0$ implies $a=0$. Let us take $a\in E$ with $A(a) =0$. It then follows that $\{A(a),A(b),A(a)\} =0=\{A(a),A(a),A(a)\}$ for all $b\in E$, and thus $\{a,b,a\} =\{a,a,a\}$ for all $b\in E$. Taking $b=0$ we get $\|a\|^{3} =\|\{a,a,a\}\| = 0$. \smallskip
	  
$(b)$ Since $A$ is a bijection and, by \Cref{l additive preservers of truncations coincide with additive preservers of inner quadratic annihilators}, $A(\leftindex^{\perp_{q}}\{a\}) = \leftindex^{\perp_{q}}\{A(a)\}$ for all $a\in E$, the desired equivalence is a consequence of the characterization of one-dimensional JB$^*$-triples given in \Cref{r chracterization of C}.\smallskip

$(c)$ This statement follows from \Cref{l additive preservers of truncations coincide with additive preservers of inner quadratic annihilators} and \Cref{th characterization of minimal tripotents}.\smallskip

$(d)$ Since $A(0) = 0,$ it suffices to prove that for each $\lambda \in \mathbb{C}\backslash \{0\}$ the element $A(\lambda e)$ lies in $\mathbb{C} A(e)$. \Cref{l additive preservers of truncations coincide with additive preservers of inner quadratic annihilators} implies that $$\leftindex^{\perp_{q}}\{A(\lambda e)\} = A ( \leftindex^{\perp_{q}}\{\lambda e\} ) = A(\leftindex^{\perp_{q}}\{e\}) = \leftindex^{\perp_{q}}\{A(e)\}.$$ Since, by $(b)$ and $(a),$ the element $A(\lambda e)$ is a positive multiple of a minimal tripotent in $F$, \Cref{l lineraly independence of multiples of minimal tripotents} assures that $A(\lambda e) \in    \mathbb{C} A(e)$. The rest follows from the bijectivity of $A$.\smallskip

$(e)$ We know from $(d)$ that $$A (e) = \gamma r(A(e)) \hbox{ and } A^{-1} (r(A(e))) = \delta r(A^{-1} (r(A(e)))),$$ for some positive real numbers $\gamma$ and $\delta$, where $r(A(e))$ and $ r(A^{-1} (r(A(e))))$ are minimal tripotents in $F$ and $E$, respectively. It follows from $(c),$ applied to $A^{-1}$, that $e = A^{-1} A(e) = A^{-1} (\gamma r(A(e))) \in \mathbb{C} A^{-1}(r(A(e))),$ which implies that $A^{-1}(r(A(e)))$ $= \lambda e$ for some non-zero complex number $\lambda$, and hence $r(A^{-1}(r(A(e))))\in  \mathbb{T} e$, as desired.
\end{proof}

\begin{rem}\label{r existence of gammaAe and gamm1Ar} Let $E$ and $F$ be JB$^*$-triples. Thanks to the previous lemma, we can conclude that every surjective additive mapping $A: E\to F$ preserving truncations in both directions is a bijection (and $A^{-1}$ enjoys the same property).  Moreover, for each minimal tripotent $e\in E$, the range tripotent of $A(e)$ (in $F^{**}$) is a minimal tripotent and lies in $F$, and there exists a unique positive number $\boldsymbol{\gamma}(e,A)$ (depending on $e$ and $A$) such that $A(e) = \boldsymbol{\gamma}(e,A) r(A(e))$. We further know the existence of a unique  $\boldsymbol{\gamma}_{1}(e,A)\in \mathbb{T}$ satisfying $A^{-1} (r(A(e))) = \boldsymbol{\gamma}_{1}(e,A)  r(A^{-1} (r(A(e))))$. We keep this notation henceforth.
\end{rem}	

A \emph{quadrangle} in $E$ is an ordered quadruple $(u_{1},u_{2},u_{3},u_{4})$ formed by tripotents in $E$ satisfying $u_{1}\bot u_{3}$, $u_{2}\bot u_{4}$, $u_{1}\top u_{2}$ $\top u_{3}\top u_{4}$ $\top u_{1}$ and $u_{4}=2 \{{u_{1}},{u_{2}},{u_{3}}\}$ (the latter equality also holds if the indices are permutated cyclically --e.g. $u_{2} = 2 \{{u_{3}},{u_{4}},{u_{1}}\}$). An ordered triplet $(w_1,u,w_2)$ of tripotents in $E$ is called a \emph{trangle} if $w_1\bot w_2$, $u\vdash w_1$, $u\vdash w_2$ and $ w_1 = Q(u) (w_2)$ (see, the references \cite{Dang_Friedm_MathScand_1987,Neher87} for additional details). \smallskip

Let $(e_1,e_2,e_3,e_4)$ (respectively, $(w_1,u,w_2)$) be a quadrangle (respectively, a trangle) of tripotents in a JB$^*$-triple $E$. It is part of the folklore in JB$^*$-triple theory that in case that $E$ is a Cartan factor and $e_1,$ $e_2$, $e_3,$ and $e_2$ are minimal tripotents (respectively, $u$ is a rank-$2$ tripotent and $w_1$ and $w_2$ are minimal tripotents), then the element $\alpha e_1 + \beta e_2 +\gamma e_4 + \delta e_3$ (respectively, $\alpha w_1 + \beta u + \delta w_2$) is a minimal tripotent in $E$ for every $\alpha,\beta,\gamma,\delta \in \mathbb{C}$ with $|\alpha|^2 + |\beta|^2 + |\gamma|^2+ |\delta|^2 =1$, and $\alpha \delta -\beta \gamma = 0$ (respectively, for every $\alpha,\beta,\delta \in \mathbb{C}$ with $|\alpha|^2 + 2 |\beta|^2 + |\delta|^2 =1$, and $\alpha \delta -\beta^2= 0$) (cf. \cite[Lemma 3.10]{Polo_Peralta_AdvMath_2018}). We can now prove a stronger property via an alternative and self-contained argument.

\begin{lem}\label{r linear combinatons of trangles and quadrangles} Let $E$ be a JB$^*$-triple. The the following statements hold:
\begin{enumerate}[$(a)$] 
	\item Let $(e_1,e_2,e_3,e_4)$ be a quadrangle of tripotents in $E$. Then for every $\alpha,\beta,\gamma,\delta \in \mathbb{C}$ with $|\alpha|^2 + |\beta|^2 + |\gamma|^2+ |\delta|^2 =1,$ and $\alpha \delta -\beta \gamma = 0$, the element $v = \alpha e_1 + \beta e_2 +\gamma e_4 + \delta e_3$ is a tripotent in $E$. Furthermore, if $e_1$ and $e_3$ {\rm(}respectively, $e_2$ and $e_4${\rm)} are minimal tripotents, the tripotents $v, e_2,$ and $e_4$ {\rm(}respectively, $v, e_1,$ and $e_3${\rm)} are minimal. 
	\item Let $(w_1,u,w_2)$ be a trangle of tripotents in $E$. Then for every $\alpha,\beta,\delta \in \mathbb{C}$ with $|\alpha|^2 + 2 |\beta|^2 + |\delta|^2 =1,$ and $\alpha \delta -\beta^2 = 0,$ the element $v= \alpha w_1 + \beta u + \delta w_2$ is a tripotent in $E$. Furthermore, if $w_1$ and $w_2$ are minimal tripotents, the tripotent $v$ is minimal. 
\end{enumerate}   	
\end{lem}

\begin{proof} $(a)$ Let $(e_1,e_2,e_3,e_4)$ be a quadrangle of tripotents in $E$, and let $\alpha,\beta,\gamma,\delta$ be complex numbers satisfying the hypotheses above. Having in mind the definition of quadrangle and Peirce arithmetic, and a bit of patience, we compute the cube of the element $v$,
\begin{align*}\{v,v,v\} &= 
\left\{ \alpha e_1 + \beta e_2 +\gamma e_4 + \delta e_3, \alpha e_1 + \beta e_2 +\gamma e_4 + \delta e_3, \alpha e_1 + \beta e_2 +\gamma e_4 + \delta e_3 \right\} \\
&= |\alpha|^2 \Big(\alpha  \{e_1,e_1, e_1\} +  \beta \{e_1,e_1,  e_2 \} + \gamma  \{e_1,e_1, e_4 \} + \delta  \CancelTo[\color{blue}]{\hspace*{-.9cm} \hbox{\scalebox{.7}{$0$  $\perp$}}}{\{e_1,e_1, e_3 \}} \Big)  
\end{align*} \vspace*{-3.9mm}
\begin{align*}
	&+ \alpha \overline{\beta} \Big(\alpha  \CancelTo[\color{blue}]{\hspace*{-1.5cm} \hbox{\scalebox{.7}{ $0$ Peirce}}}{\{e_1,e_2, e_1\}} +  \beta \{e_1,e_2,  e_2 \} + \gamma  \CancelTo[\color{blue}]{0}{\{e_1,e_2, e_4 \}} + \delta  \{e_1,e_2, e_3 \} \Big) 
	\end{align*} \vspace*{-3.9mm}
	\begin{align*}
	&+ \alpha \overline{\gamma} \Big(\alpha   \CancelTo[\color{blue}]{\hspace*{-1.5cm} \hbox{\scalebox{.7}{ $0$ Peirce}}}{\{e_1,e_4, e_1\}} +  \beta \CancelTo[\color{blue}]{\hspace*{-.9cm} \hbox{\scalebox{.7}{$0$  $\perp$}}}{\{e_1,e_4,  e_2 \}} + \gamma  \{e_1,e_4, e_4 \} + \delta  \{e_1,e_4, e_3 \} \Big) \\
&+ \alpha \overline{\delta} \Big(\alpha  \CancelTo[\color{blue}]{\hspace*{-.9cm} \hbox{\scalebox{.7}{$0$  $\perp$}}}{\{e_1,e_3, e_1\}} +  \beta \CancelTo[\color{blue}]{\hspace*{-.9cm} \hbox{\scalebox{.7}{$0$  $\perp$}}}{\{e_1,e_3,  e_2 \}} + \gamma  \CancelTo[\color{blue}]{\hspace*{-.9cm} \hbox{\scalebox{.7}{$0$  $\perp$}}}{\{e_1,e_3, e_4 \}} + \delta  \CancelTo[\color{blue}]{\hspace*{-.9cm} \hbox{\scalebox{.7}{$0$  $\perp$}}}{\{e_1,e_3, e_3 \}} \Big)
\end{align*} \vspace*{-3.8mm}
\begin{align*}
	&+ \beta \overline{\alpha} \Big(\alpha  \{e_2,e_1, e_1\} +  \beta \CancelTo[\color{blue}]{\hspace*{-1.5cm} \hbox{\scalebox{.7}{ $0$ Peirce}}}{\{e_2,e_1,  e_2 \}} + \gamma  \{e_2,e_1, e_4 \} + \delta  \CancelTo[\color{blue}]{\hspace*{-.9cm} \hbox{\scalebox{.7}{$0$  $\perp$}}}{\{e_2,e_1, e_3 \}} \Big)\\
	&+  \beta \overline{\beta} \Big(\alpha  \{e_2,e_2, e_1\} +  \beta \{e_2,e_2,  e_2 \} + \gamma  \CancelTo[\color{blue}]{\hspace*{-.9cm} \hbox{\scalebox{.7}{$0$  $\perp$}}}{\{e_2,e_2, e_4 \}} + \delta  \{e_2,e_2, e_3 \} \Big) 
	\end{align*} \vspace*{-3.8mm}
	\begin{align*}
	&+  \beta \overline{\gamma} \Big(\alpha  \CancelTo[\color{blue}]{\hspace*{-.9cm} \hbox{\scalebox{.7}{$0$  $\perp$}}}{\{e_2,e_4, e_1\}} +  \beta \CancelTo[\color{blue}]{\hspace*{-.9cm} \hbox{\scalebox{.7}{$0$  $\perp$}}}{\{e_2,e_4,  e_2 \}} + \gamma  \CancelTo[\color{blue}]{\hspace*{-.9cm} \hbox{\scalebox{.7}{$0$  $\perp$}}}{\{e_2,e_4, e_4 \}} + \delta  \CancelTo[\color{blue}]{\hspace*{-.9cm} \hbox{\scalebox{.7}{$0$  $\perp$}}}{\{e_2,e_4, e_3 \}} \Big)\\
	&+  \beta \overline{\delta} \Big(\alpha  \CancelTo[\color{blue}]{\hspace*{-.9cm} \hbox{\scalebox{.7}{$0$  $\perp$}}}{\{e_2,e_3, e_1\}} +  \beta \CancelTo[\color{blue}]{\hspace*{-1.5cm} \hbox{\scalebox{.7}{ $0$ Peirce}}}{\{e_2,e_3,  e_2 \}} + \gamma  \{e_2,e_3, e_4 \} + \delta  \{e_2,e_3, e_3 \} \Big) 
\end{align*} \vspace*{-3.8mm}
\begin{align*}
&+ \gamma \overline{\alpha} \Big(\alpha  \{e_4,e_1, e_1\} +  \beta \{e_4,e_1,  e_2 \} + \gamma  \CancelTo[\color{blue}]{\hspace*{-1.5cm} \hbox{\scalebox{.7}{ $0$ Peirce}}}{\{e_4,e_1, e_4 \}} + \delta  \CancelTo[\color{blue}]{\hspace*{-.9cm} \hbox{\scalebox{.7}{$0$  $\perp$}}}{\{e_4,e_1, e_3 \}} \Big) \\ 
&+  \gamma \overline{\beta} \Big(\alpha  \CancelTo[\color{blue}]{\hspace*{-.9cm} \hbox{\scalebox{.7}{$0$  $\perp$}}}{\{e_4,e_2, e_1\}} +  \beta \CancelTo[\color{blue}]{\hspace*{-.9cm} \hbox{\scalebox{.7}{$0$  $\perp$}}}{\{e_4,e_2,  e_2 \}} + \gamma  \CancelTo[\color{blue}]{\hspace*{-.9cm} \hbox{\scalebox{.7}{$0$  $\perp$}}}{\{e_4,e_2, e_4 \}} + \delta  \CancelTo[\color{blue}]{\hspace*{-.9cm} \hbox{\scalebox{.7}{$0$  $\perp$}}}{\{e_4,e_2, e_3 \}} \Big) 
\end{align*} \vspace*{-3.8mm}
\begin{align*} 
&+  \gamma \overline{\gamma} \Big(\alpha  \{e_4,e_4, e_1\} +  \beta \CancelTo[\color{blue}]{\hspace*{-.9cm} \hbox{\scalebox{.7}{$0$  $\perp$}}}{\{e_4,e_4,  e_2 \}} + \gamma  \{e_4,e_4, e_4 \} + \delta  \{e_4,e_4, e_3 \} \Big) \\ 
&+ \gamma \overline{\delta} \Big(\alpha  \CancelTo[\color{blue}]{\hspace*{-.9cm} \hbox{\scalebox{.7}{$0$  $\perp$}}}{\{e_4,e_3, e_1\}} +  \beta \{e_4,e_3,  e_2 \} + \gamma  \CancelTo[\color{blue}]{\hspace*{-1.5cm} \hbox{\scalebox{.7}{ $0$ Peirce}}}{\{e_4,e_3, e_4 \}} + \delta  \{e_4,e_3, e_3 \} \Big) 
\end{align*} \vspace*{-4.1mm}
\begin{align*}
	&+ \delta \overline{\alpha} \Big(\alpha  \CancelTo[\color{blue}]{\hspace*{-.9cm} \hbox{\scalebox{.7}{$0$  $\perp$}}}{\{e_3,e_1, e_1\}} +  \beta \CancelTo[\color{blue}]{\hspace*{-.9cm} \hbox{\scalebox{.7}{$0$  $\perp$}}}{\{e_3,e_1,  e_2 \}} + \gamma  \CancelTo[\color{blue}]{\hspace*{-.9cm} \hbox{\scalebox{.7}{$0$  $\perp$}}}{\{e_3,e_1, e_4 \}} + \delta  \CancelTo[\color{blue}]{\hspace*{-.9cm} \hbox{\scalebox{.7}{$0$  $\perp$}}}{\{e_3,e_1, e_3 \}} \Big) \\ 
	&+  \delta \overline{\beta} \Big(\alpha  \{e_3,e_2, e_1\} +  \beta \{e_3,e_2,  e_2 \} + \gamma  \CancelTo[\color{blue}]{\hspace*{-.9cm} \hbox{\scalebox{.7}{$0$  $\perp$}}}{\{e_3,e_2, e_4 \}} + \delta  \CancelTo[\color{blue}]{\hspace*{-1.5cm} \hbox{\scalebox{.7}{ $0$ Peirce}}}{\{e_3,e_2, e_3 \}} \Big) 
	\end{align*} \vspace*{-4.3mm}
	\begin{align*}
	&+  \delta \overline{\gamma} \Big(\alpha  \{e_3,e_4, e_1\} +  \beta \CancelTo[\color{blue}]{\hspace*{-.9cm} \hbox{\scalebox{.7}{$0$  $\perp$}}}{\{e_3,e_4,  e_2 \}} + \gamma  \{e_3,e_4, e_4 \} + \delta  \CancelTo[\color{blue}]{\hspace*{-1.5cm} \hbox{\scalebox{.7}{ $0$ Peirce}}}{\{e_3,e_4, e_3 \}} \Big) \\ 
	&+ \delta \overline{\delta} \Big(\alpha  \CancelTo[\color{blue}]{\hspace*{-.9cm} \hbox{\scalebox{.7}{$0$  $\perp$}}}{\{e_3,e_3, e_1\}} +  \beta \{e_3,e_3,  e_2 \} + \gamma  \{e_3,e_3, e_4 \} + \delta  \{e_3,e_3, e_3 \} \Big)
\end{align*} \vspace*{-3.8mm}
\begin{align*}
=&  |\alpha|^2 \alpha  e_1+ \frac12 |\alpha|^2 \beta e_2  + \frac12 |\alpha|^2 \gamma   e_4  +\frac12 \alpha \overline{\beta}  \beta e_1 + \frac12 \alpha \overline{\beta} \delta  e_4 + \frac12 \alpha \overline{\gamma}  \gamma  e_1 + \frac12 \alpha \overline{\gamma} \delta  e_2 \\
& +\frac12 \beta \overline{\alpha} \alpha  e_2 +  \frac12 \beta \overline{\alpha} \gamma  e_3 + \frac12 \beta \overline{\beta} \alpha  e_1 + \beta \overline{\beta} \beta e_2 + \frac12 \beta \overline{\beta} \delta  e_3 + \frac12 \beta \overline{\delta} \gamma  e_1 +\frac12 \beta \overline{\delta} \delta  e_2  \\ 
&+ \frac12 \gamma \overline{\alpha} \alpha  e_4 +\frac12 \gamma \overline{\alpha} \beta e_3 + \frac12 \gamma \overline{\gamma} \alpha  e_1 + \gamma \overline{\gamma} \gamma  e_4 +  \frac12 \gamma \overline{\gamma} \delta  e_3 + \frac12 \gamma \overline{\delta} \beta e_1 + \frac12 \gamma \overline{\delta} \delta  e_4  \\
& + \frac12 \delta \overline{\beta} \alpha  e_4 + \frac12 \delta \overline{\beta} \beta e_3+  \frac12 \delta \overline{\gamma} \alpha  e_2 + \frac12 \delta \overline{\gamma} \gamma e_3 + \frac12 \delta \overline{\delta}  \beta e_2  +\frac12  \delta \overline{\delta} \gamma  e_4 + \delta \overline{\delta} \delta  e_3  
\end{align*} \vspace*{-3.9mm}
\begin{align*}
&= \left(|\alpha|^2 \alpha  + \alpha |\beta|^2 + \alpha |\gamma|^2  + \beta \overline{\delta} \gamma  \right) e_1 + \left( |\alpha|^2 \beta + \alpha \overline{\gamma} \delta  + \beta |\beta|^2 +\beta |\delta|^2   \right) e_2 \hspace*{9mm} \\
&+ \left( \beta \overline{\alpha} \gamma  + |\beta|^2 \delta  + |\gamma|^2  \delta   + \delta |\delta|^2  \right)   e_3 + \left(|\alpha|^2 \gamma  +\alpha \overline{\beta} \delta + \gamma |\gamma|^2 + \gamma |\delta|^2\right)   e_4 \\
&= \alpha e_1 + \beta e_2 +\gamma e_4 + \delta e_3 =v.
\end{align*}	

Since, clearly, $e_1, e_2, e_4, e_3\in E_2(e_1+ e_3)$, all the tripotents listed above lie in $E_2(e_1+ e_3)$.  If $e_1$ and $e_3$ are minimal tripotents, the Peirce-$2$ subspace $E_2(e_1+e_3)$ is a spin factor (cf. \cite[Lemma 3.6]{Kal_Peralta_Ann_Math_Phys_2021}). Since spin factors have rank-$2$ (see, for example, \cite[Table~1 in page 210]{Kaup_ManMath_1997}), $e_2\perp e_4,$ and $e_2+e_4\in E_2 (e_1+e_3)$, it follows that $e_2$ and $e_4$ must be minimal in $E_2(e_1+e_2)$, and hence in $E$ by \cite[Lemma~3.3]{CabMarPe2024}. Consequently, $e_2$ and $e_4$ are minimal tripotents too. It is part of the folklore in JB$^*$-triple theory that under these conditions $v$ is a minimal tripotent.\smallskip

We can alternatively consider the element $\tilde{v} = \overline{\delta} e_1 -\overline{\gamma}  e_2 - \overline{\beta} e_4 +  \overline{\alpha}  e_3$. The arguments above also prove that $\tilde{v}$ is a tripotent in $E$. By following similar computations to those at the beginning of this proof we deduce that \begin{align*}
\{v,v,\tilde{v}\} &= \{\alpha e_1 + \beta e_2 +\gamma e_4 + \delta e_3, \alpha e_1 + \beta e_2 +\gamma e_4 + \delta e_3, \overline{\delta} e_1 -\overline{\gamma}  e_2 - \overline{\beta} e_4 +  \overline{\alpha}  e_3\} \\
=&  |\alpha|^2 \overline{\delta}  e_1- \frac12 |\alpha|^2 \overline{\gamma} e_2  - \frac12 |\alpha|^2 \overline{\beta}   e_4  -\frac12 \alpha \overline{\beta}  \overline{\gamma} e_1 + \frac12 \alpha \overline{\beta} \overline{\alpha}  e_4 - \frac12 \alpha \overline{\gamma} \overline{\beta}  e_1 + \frac12 \alpha \overline{\gamma} \overline{\alpha}  e_2 \\
& +\frac12 \beta \overline{\alpha} \overline{\delta}  e_2 -  \frac12 \beta \overline{\alpha} \overline{\beta}  e_3 + \frac12 \beta \overline{\beta} \overline{\delta}  e_1 - \beta \overline{\beta} \overline{\gamma} e_2 + \frac12 \beta \overline{\beta} \overline{\alpha}  e_3 - \frac12 \beta \overline{\delta} \overline{\beta}  e_1 +\frac12 \beta \overline{\delta} \overline{\alpha}  e_2 
\end{align*}\vspace*{-3.9mm}
\begin{align*}
	\hspace*{9mm}& + \frac12 \gamma \overline{\alpha} \overline{\delta}  e_4 -\frac12 \gamma \overline{\alpha} \overline{\gamma} e_3 + \frac12 \gamma \overline{\gamma} \overline{\delta}  e_1 - \gamma \overline{\gamma} \overline{\beta}  e_4 +  \frac12 \gamma \overline{\gamma} \overline{\alpha}  e_3 - \frac12 \gamma \overline{\delta} \overline{\gamma} e_1 + \frac12 \gamma \overline{\delta} \overline{\alpha}  e_4  \\
	& + \frac12 \delta \overline{\beta} \overline{\delta} e_4 - \frac12 \delta \overline{\beta} \overline{\gamma} e_3+  \frac12 \delta \overline{\gamma} \overline{\delta}  e_2 - \frac12 \delta \overline{\gamma} \overline{\beta} e_3 - \frac12 \delta \overline{\delta}  \overline{\gamma} e_2  -\frac12  \delta \overline{\delta} \overline{\beta}  e_4 + \delta \overline{\delta} \overline{\alpha}  e_3 \\
\end{align*}\vspace*{-9.9mm}
\begin{align*} &= \left(  |\alpha|^2 \overline{\delta} -\frac12 \alpha \overline{\beta}  \overline{\gamma} - \frac12 \alpha \overline{\gamma} \overline{\beta} + \frac12 \beta \overline{\beta} \overline{\delta}  - \frac12 \beta \overline{\delta} \overline{\beta} + \frac12 \gamma \overline{\gamma} \overline{\delta} - \frac12 \gamma \overline{\delta} \overline{\gamma}  \right)  e_1 \hspace{12mm}\\
	&\hspace{2mm} + \left( - \frac12 |\alpha|^2 \overline{\gamma} + \frac12 \alpha \overline{\gamma} \overline{\alpha} +\frac12 \beta \overline{\alpha} \overline{\delta}  - \beta \overline{\beta} \overline{\gamma} +\frac12 \beta \overline{\delta} \overline{\alpha} +  \frac12 \delta \overline{\gamma} \overline{\delta} - \frac12 \delta \overline{\delta}  \overline{\gamma}  \right)  e_2 \\
	&\hspace{2mm} + \left(  -  \frac12 \beta \overline{\alpha} \overline{\beta}   + \frac12 \beta \overline{\beta} \overline{\alpha}  -\frac12 \gamma \overline{\alpha} \overline{\gamma} +  \frac12 \gamma \overline{\gamma} \overline{\alpha} - \frac12 \delta \overline{\beta} \overline{\gamma}  - \frac12 \delta \overline{\gamma} \overline{\beta}  + \delta \overline{\delta} \overline{\alpha} \right) e_3  \\
	&\hspace{2mm} + \left(- \frac12 |\alpha|^2 \overline{\beta} + \frac12 \alpha \overline{\beta} \overline{\alpha} + \frac12 \gamma \overline{\alpha} \overline{\delta}  - \gamma \overline{\gamma} \overline{\beta}  + \frac12 \gamma \overline{\delta} \overline{\alpha} + \frac12 \delta \overline{\beta} \overline{\delta} -\frac12  \delta \overline{\delta} \overline{\beta}   \right) e_4 =0,
\end{align*}  and thus $v\perp \tilde{v}$. Since $v,\tilde{v}$ are both non-zero and lie in $E_2(e_1+e_3)$, we conclude that $v$ and $\tilde{v}$ are minimal tripotents in $E$.\smallskip 

Clearly, the roles of $e_1,e_3$ and $e_2,e_4$ can be exchanged in this final argument by just having in mind that $(e_2,e_1,e_4,e_3)$ is a quadrangle too.\smallskip

$(b)$ The proof follows similar lines to those in the proof of $(a)$ by just observing the following identities. First, since by assumptions $\{u,w_1,u\} = w_2$, the element $w_1+w_2$ is a self-adjoint tripotent in the JB$^*$-algebra $E_2(u)$. Moreover, since $\{w_1+w_2,w_1+w_2,u\} =  \{w_1,w_1,u\} + \{w_2,w_2,u\} = u$, we can deduce, via Peirce arithmetic, that $u = \{w_1+w_2,u, w_1+w_2\} = \CancelTo[\color{blue}]{\hspace*{-.5cm} \hbox{\scalebox{.7}{ $0$}}}{\{w_1,u, w_1\}}+ \CancelTo[\color{blue}]{\hspace*{-.5cm} \hbox{\scalebox{.7}{ $0$}}}{\{w_2,u,w_2\}} + 2 \{w_1,u, w_2\} = 2 \{w_1,u, w_2\}.$ For the sake of brevity, the remaining computations are left to the reader. 
\end{proof}

The relationship ``being in the inner quadratic annihilator'' is not, in general, reflexive. For example, the tripotents $e = \left(\begin{matrix} 1 & 0 \\
	0 & 0 
\end{matrix}\right)$ and $w = \left(\begin{matrix} 0 & 1 \\
1 & 0 
\end{matrix}\right)$ satisfy that $w\in  \leftindex^{\perp_{q}}\{e\}$ but $e\notin \leftindex^{\perp_{q}}\{w\}$.\smallskip

It is shown in \cite[Lemma in page 306]{Dang_Friedm_MathScand_1987} that each family $\{w_j\}_j$ of mutually collinear minimal tripotents in a JB$^*$-triple generates a norm-closed subspace isometrically isomorphic to a complex Hilbert space in which the family $\{w_j\}_j$ is an orthonormal basis. We shall see a kind of reciprocal to this property in the next key proposition. 

\begin{prop}\label{prop characterization of collinearity} Let $w_1,w_2$ be two minimal tripotents in a JB$^*$-triple $E$. Then the following statements hold:
	\begin{enumerate}[$(a)$]
		\item $w_1\in \leftindex^{\perp_{q}}\{w_2\}$ if, and only if, $w_2\in \leftindex^{\perp_{q}}\{w_1\}.$
		\item Suppose that $w_2\in \leftindex^{\perp_{q}}\{w_1\}$ and there exist $\gamma_1,\gamma_2\in \mathbb{R}^+,$ $q_1,q_2\in \mathbb{Q}$ with $q_1^2 + q_2^2 =1,$ and $q_1 q_2\neq 0$, such that the element $q_1 \gamma_1 w_1 + q_2 \gamma_2 w_2$ is a positive multiple of a minimal tripotent in $E$. Then $w_1$ and $w_2$ are collinear. 
	\end{enumerate}
\end{prop}

\begin{proof} As in the proof of \Cref{th characterization of minimal tripotents}, we can embed $E$ inside the atomic part, $\mathcal{A}_{E^{**}},$ of its bidual space, $E^{**}$, via an isometric triple isomorphism $\Psi_{E}: E \hookrightarrow \mathcal{A}_{E^{**}}$ with weak$^*$-dense image, and $\mathcal{A}_{E^{**}}$ can be represented as an $\ell_{\infty}$-sum of a family of Cartan factors $\{C_i\}_{i}$ (cf. \cite[Theorem 2]{FriedRusso_Predual} and \cite[Proposition 1 and its proof]{FriedRusso_GelfNaim}).  Clearly, $\Psi_{E} (w_1)$ and $\Psi_{E} (w_2)$ are minimal tripotents in $\mathcal{A}_{E^{**}} = \bigoplus^{\infty} C_i$, and hence each one of them belongs to a unique Cartan factor in the $\ell_{\infty}$-sum. Let us assume that $\Psi_{E} (w_j) \in C_{i_{j}}$ ($j=1,2$). 
	\smallskip  
	
	$(a)$ Assume that $w_2\in \leftindex^{\perp_{q}}\{w_1\}$ (equivalently, $\Psi_{E}(w_2)\in \leftindex^{\perp_{q}}\{\Psi_{E}(w_1)\}$). If $i_1\neq i_2$, the tripotents $\Psi_{E}(w_1)$ and $\Psi_{E}(w_2)$ (equivalently, $w_1,w_2$) are orthogonal in $\mathcal{A}_{E^{**}}$ (respectively, in $E$), which implies that $w_1\in \leftindex^{\perp_{q}}\{w_2\}$ and $w_2\in \leftindex^{\perp_{q}}\{w_1\}.$ We can therefore assume that $\Psi_{E}(w_1)$ and $\Psi_{E}(w_2)$ belong to the same Cartan factor ${C}_{i_1}$ (i.e. $i_1 = i_2$). We are thus, in a position to apply \cite[Lemma 3.10]{Polo_Peralta_AdvMath_2018} to deduce that one of the following statements holds:
	
\noindent $(i)$ There exist minimal tripotents $e_2,e_3,e_4$ in ${C}_{i_1}$ such that $(\Psi_{E}(w_1),e_2,e_3,e_4)$ is a quadrangle and $\Psi_{E}(w_2) = \alpha \Psi_{E}(w_1) + \beta e_2 + \gamma e_4 + \delta e_3$ with $\alpha, \beta, \gamma, \delta\in \mathbb{C},$ $|\alpha|^2 + |\beta|^2 +|\gamma|^2+ |\delta|^2 =1,$ and $\alpha \delta = \gamma \beta$;\smallskip

\noindent $(ii)$ There exist a minimal tripotent $\tilde e\in {C}_{i_1}$ and a rank-two tripotent $u\in {C}_{i_1}$ such that $(\Psi_{E}(w_1), u,\tilde e)$ is a trangle and $\Psi_{E}(w_2) = \alpha \Psi_{E}(w_1) + \beta u + \delta \tilde{e}$ with $\alpha, \beta, \delta\in \mathbb{C}$ and $\alpha \delta = \beta^2$.
\smallskip 
	
	Since $\Psi_{E}(w_2)\in \leftindex^{\perp_{q}}\{\Psi_{E}(w_1)\}$, it follows that in each of the previous cases we have $\alpha =0$. In case $(ii)$ the element $\beta$ must be also zero, and hence $\Psi_{E}(w_2) = \delta \tilde{e}$ is orthogonal to $\Psi_{E}(w_1),$ and consequently, $\Psi_{E}(w_1)\in  \leftindex^{\perp_{q}}\{\Psi_{E}(w_2)\}$. In case $(i)$, $\beta \gamma =0$, implies that $\Psi_{E}(w_2) =  \gamma e_4 + \delta e_3$ or $\Psi_{E}(w_2) =  \beta e_2 + \delta e_3$. In both cases, by Peirce arithmetic, we have $\{ \Psi_{E}(w_2), \Psi_{E}(w_1), \Psi_{E}(w_2) \} = \{ \gamma e_4 + \delta e_3, \Psi_{E}(w_1), \gamma e_4 + \delta e_3 \} = \gamma \{  e_4 , \Psi_{E}(w_1), e_4 \} = 0$, and similarly $\{ \beta e_2 + \delta e_3, \Psi_{E}(w_1), \beta e_2 + \delta e_3 \} = \beta^2 \{ e_2 , \Psi_{E}(w_1), e_2 \}=0$. This shows that $\Psi_{E}(w_1)\in \leftindex^{\perp_{q}}\{\Psi_{E}(w_2)\}$. \smallskip

	$(b)$ Let us assume the hypothesis in the statement and in the opening paragraph of this proof. If $i_1\neq i_2$, the tripotents $w_1$ and $w_2$ (equivalently, $\Psi_{E}(w_1) =w_1$ and $\Psi_{E} (w_2) =w_1$) are orthogonal in $\mathcal{A}_{E^{**}}$ (respectively, in $E$). By assumptions, the element $q_1 \gamma_1 w_1 + q_2\gamma_2 w_2$ is a minimal tripotent in $E$, which is impossible. Therefore the elements $\Psi_{E}(w_1) =w_1$ and $\Psi_{E}(w_2)=w_2$ belong to the same Cartan factor ${C}_{i_1}$. Lemma 3.10 in \cite{Polo_Peralta_AdvMath_2018} assures that one of the following statements holds:\medskip
	
\noindent $(i)$ There exist minimal tripotents $e_2,e_3,e_4$ in ${C}_{i_1}$ such that the elements $w_1,$ $e_2,$ $e_3,$ and $e_4$ form a quadrangle, and $$w_2 = \alpha w_1 + \beta e_2 + \gamma e_4 + \delta e_3$$ with $\alpha, \beta, \gamma, \delta\in \mathbb{C},$ $|\alpha|^2 + |\beta|^2 +|\gamma|^2+ |\delta|^2 =1,$ and $\alpha \delta = \gamma \beta$. In such a case we have $$\begin{aligned}
			q_1 \gamma_1 w_1 + q_2\gamma_2 w_2 &=  q_1 \gamma_1 w_1 + q_2 \gamma_2 \alpha w_1 + q_2 \gamma_2 \beta e_2 + q_2 \gamma_2 \gamma e_4 + q_2 \gamma_2 \delta e_3.
		\end{aligned}$$ Now, by the quadrangle's properties, we get $\left( q_1 \gamma_1 + q_2 \gamma_2 \alpha\right) q_2 \gamma_2 \delta =  q_2^2 \gamma_2^2  \beta  \gamma, $ and thus $\gamma_1 \gamma_2 q_1 q_2 \delta =0$, equivalently, $\delta =0$. The condition $w_2\in \leftindex^{\perp_{q}}\{w_1\}$ implies that $\alpha =0$.  Therefore, either $w_2 =  \beta e_2,$ ($\gamma =0$) or $w_2 = \gamma e_4$ ($\beta =0$), and hence $w_1\top w_2$.\medskip
		
\noindent $(ii)$ There exist a minimal tripotent $\tilde e\in {C}_{i_1}$ and a rank-two tripotent $u\in {C}_{i_1}$ such that $(\Psi_{E}(w_1), u,\tilde e)$ is a trangle and $\Psi_{E}(w_2) = \alpha \Psi_{E}(w_1) + \beta u + \delta \tilde{e}$ with $\alpha, \beta, \delta\in \mathbb{C}$ and $\alpha \delta = \beta^2$. The condition $w_2\in \leftindex^{\perp_{q}}\{w_1\}$ gives $\alpha =0 = \beta,$ and thus $w_2 = \delta \tilde{e}$ is orthogonal to $w_1$, which contradicts that $q_1 \gamma_1 w_1 + q_2 \gamma_2 w_2$ is a positive scalar multiple of a minimal tripotent in $E$. This shows that this second case is impossible. 
\end{proof}

As a consequence of the previous proposition, we can conclude next that a surjective additive mapping between JB$^*$-triples preserving truncations in both directions ``somehow preserves'' collinear minimal tripotents.

\begin{prop}\label{prop preservation of collinearity} Let $A: E\to F$ be a surjective additive mapping between JB$^*$-triples. Suppose additionally that $A$ preserves truncations of elements in both directions. Let $w_1,w_2$ be two collinear minimal tripotents in $E$. Then the range tripotents $r(A(w_1))$ and $r(A(w_2))$ are minimal and collinear in $F$.
\end{prop}

\begin{proof} Since $w_1$ and $w_2$ are collinear, it is not hard to see that for every $q_1,q_2\in \mathbb{Q}$ with $q_1^2 + q_2^2 =1$, the element $q_1 w_1 + q_2 w_2$ is a minimal tripotent in $E$, and thus, \Cref{r existence of gammaAe and gamm1Ar} assures that $$A(q_1 w_1 + q_2 w_2) = q_1 \boldsymbol{\gamma}(w_1,A) r(A(w_1)) + q_2 \boldsymbol{\gamma}(w_2, A) r(A(w_2))$$ is a positive scalar multiple of a minimal tripotent in $F$. It follows from \Cref{prop characterization of collinearity}$(b)$ that $r(A(w_1)) \top r(A(w_2))$.	
\end{proof}

We say that two tripotents $e$ and $v$ in a JB$^*$-triple $E$ are \emph{compatible} if $P_j (e) P_k(v) = P_k (v) P_j (e)$ for all $i,k\in\{0,1,2\}$. If $v\in E_j(e)$ for some $j\in\{0,1,2\}$, then $e$ and $v$ are compatible \cite[$(1.10)$]{Horn_MathScand_1987}. 

\begin{lem}\label{l orthogonal tripotents do not admit a collinear to one of them but not to the other} Let $w_1,w_2$ be orthogonal tripotents in a JB$^*$-triple $E$. Then $$E_1 (w_1) \subseteq   \leftindex^{\perp_{q}}\{w_2\}.$$  
\end{lem}

\begin{proof} Let us take $x\in E_1 (w_1)$. Since $w_1$ and $w_2$ are compatible, we conclude that $E_1 (w_1) = \left(E_1 (w_1)\cap E_0 (w_2)\right) \oplus \left(E_1 (w_1)\cap E_1 (w_2)\right) \oplus \left(E_1 (w_1)\cap E_2 (w_2)\right).$ Observe that $w_2\perp w_1$ implies that $E_2 (w_2) \subseteq E_0 (w_1)$, and thus $E_1 (w_1)\cap E_2 (w_2) = \{0\}$. We can therefore write $x =x_0 + x_1$ with $x_j \in E_1 (w_1)\cap E_j (w_2)$. Therefore, we deduce from Peirce arithmetic that $$\{w_2, x, w_2\} = \{w_2, x_0, w_2\} + \{w_2, x_1, w_2\} = \{w_2, x_1, w_2\}\in E_{4-1} (w_2) =\{0\},$$ which shows that $x\in \leftindex^{\perp_{q}}\{w_2\}$.
\end{proof}

We can now describe the images of two orthogonal minimal tripotents under a surjective additive mapping preserving truncations. 

\begin{prop}\label{prop preservation of orthogonality} Let $E$ and $F$ be JB$^*$-triples, and let $A: E\to F$ be a surjective additive mapping preserving truncations in both directions. Suppose additionally that $F$ is atomic. Let $w_1,w_2$ be two orthogonal minimal tripotents in $E$. Then the range tripotents $r(A(w_1))$ and $r(A(w_2))$ are orthogonal in $F$.
\end{prop}

\begin{proof} Since $F$ is atomic, it can be expressed in the form $F=\bigoplus^{\ell_{\infty}} C_i,$ where each $C_i$ is a Cartan factor. If $A(w_1)$ and $ A(w_2)$ lie in two different Cartan factors in the decomposition of $F$, then they are clearly orthogonal and the same occurs to $ r(A(w_1))$ and $r(A(w_2))$, which gives the desired statement.  We can therefore assume that $A(w_1)$ and $A(w_2)$ belong to the same Cartan factor $C_{i_1}$. We are again in a position to apply \cite[Lemma 3.10]{Polo_Peralta_AdvMath_2018} to reduce our study to one of the following cases:\medskip

\noindent $(i)$ There exist minimal tripotents $e_2,e_3,e_4$ in ${C}_{i_1}$ such that the ordered quadruple $\{ r(A(w_1)),e_2, e_3, e_4\}$ is a quadrangle, and $$r(A(w_2)) = \alpha r(A(w_1)) + \beta e_2 + \gamma e_4 + \delta e_3$$ with $\alpha, \beta, \gamma, \delta\in \mathbb{C},$ $|\alpha|^2 + |\beta|^2 +|\gamma|^2+ |\delta|^2 =1,$ and $\alpha \delta = \gamma \beta$. By hypotheses, $A(w_2) = \boldsymbol{\gamma}(w_2,A) r(A(w_2))\in \leftindex^{\perp_{q}}\{A(w_1)\}$, and thus $\alpha =0$, which also implies that $\beta =0$ or $\gamma =0$. We shall only consider the first case, since the other one can be treated via similar arguments. Assume then that  $r(A(w_2)) = \gamma e_4 + \delta e_3.$ Suppose $\gamma\neq 0$. The element $e_4$ in $C_{i_1}\subset F$ satisfies $e_4\top r(A(w_1))$ and $\{r(A(w_2)), e_4, r(A(w_2)) \} = \gamma^2 e_4 + \delta\gamma e_3 \neq 0$.\smallskip
		
\Cref{prop preservation of collinearity} implies that $r(A^{-1} (e_4)) \top r(A^{-1} (r(A(w_1))))$. Observe that by \Cref{l basic properties of additive preservers of truncations into abd images of min trip}$(d)$ we have $ r(A^{-1} (r(A(w_1))) ) = e^{i \theta} w_1$ for some real $\theta$.  Now, it follows from $r(A^{-1} (e_4)) \top r(A^{-1} (r(A(w_1))) )$ that $w_1 \top r(A^{-1} (e_4))$.\smallskip

On the other hand, 	since $w_1\perp w_2$ and $A^{-1} (e_4) = \boldsymbol{\gamma}(e_4,A^{-1}) r(A^{-1} (e_4))\in E_1 (w_1),$ we deduce from \Cref{l orthogonal tripotents do not admit a collinear to one of them but not to the other} that $A^{-1} (e_4)\in \leftindex^{\perp_{q}}\{w_2\}.$ Therefore the hypotheses on the mapping $A$ lead to $e_4\in \leftindex^{\perp_{q}}\{A(w_2)\}$, which contradicts that $$\{r(A(w_2)), v_4, r(A(w_2)) \} \neq 0.$$ This implies that $\gamma =0$, and hence $$A(w_2) = \boldsymbol{\gamma}(w_2,A) \delta e_3 \perp  A(w_1) = \boldsymbol{\gamma}(w_1,A) r(A(w_1)),$$ as desired.\smallskip

\noindent $(ii)$ There exist a minimal tripotent $\tilde e\in {C}_{i_1}$ and a rank-two tripotent $u\in {C}_{i_1}$ such that $(r(A(w_1)), u,\tilde e)$ is a trangle and $r(A(w_2)) = \alpha r(A(w_1)) + \beta u + \delta \tilde{e}$ with $\alpha, \beta, \delta\in \mathbb{C}$ and $\alpha \delta = \beta^2$. The condition $w_2\in \leftindex^{\perp_{q}}\{w_1\}$ (equivalently, $A(w_2)\in \leftindex^{\perp_{q}}\{A(w_1)\}$) gives $\alpha =0 = \beta,$ and thus $A(w_2) = \boldsymbol{\gamma}(w_2,A) \delta \tilde{e}$ is orthogonal to $A(w_1) =  \boldsymbol{\gamma}(w_1,A) r(A(w_1))$.
\end{proof}

\section{Atomic JBW$^*$-triples}\label{sec: atomic}

Let $E$ be a JBW$^*$-triple with predual $E_*$. The extreme points of the closed unite ball of $E_*$ are called \emph{atoms} or \emph{pure atoms}  of $E$. The symbol $\partial_e (\mathcal{B}_{E_{*}})$ will denote the set of all atoms of $E$. A celebrated result by Y. Friedman and B. Russo shows that for each $\varphi\in \partial_e (\mathcal{B}_{E_*})$ there exists a unique minimal tripotent $e$ in $E$ satisfying $P_2 (e) (x) = \varphi (x) e$ for all $x\in E$ (cf. \cite[Proposition 4]{FriedRusso_Predual}). In case that $E$ is an atomic JBW$^*$-triple, $E_*$ is precisely the norm closure of the linear span of the atoms of $E$ \cite[Theorem 1]{FriedRusso_Predual}. In particular,  $\partial_e (\mathcal{B}_{E_*})$, equivalently $\{P_2 (e) : e\hbox{ minimal tripotent in } E \}$, is a norming set for $E$. \smallskip

A (closed) subtriple $I$ of a JB$^*$-triple $E$ is said to be an \emph{ideal} (respectively, an \emph{inner ideal}) of $E$ if $\{E,E,I\} + \{E,I,E\}\subseteq I$ (respectively, $\{I,E,I\} \subseteq I$).  If $e$ is a tripotent of $E$, it follows from Peirce arithmetic that $E_0(e)$ and $E_2(e)$ are inner ideals of $E$. A JBW$^*$-triple $E$ is called a factor if there does not exist a decomposition of $E$ as a direct sum of two non-zero ideals $I$, $J$, or equivalently, if $\{0\}$ and $E$ are the only weak$^*$-closed ideals of $E$. The Cartan factors are precisely the JBW$^*$-triple factors $E$ containing a minimal tripotent (cf. \cite[Corollary 1.8]{Horn_MathZ_1987}). It is known that every tripotent $e$ in a Cartan factor $C$ admits a representation $e= \sum_{i} e_i$, where the series converges with respect to the weak$^*$-topology of $C$  and $\{e_i\}_{i}$ is a family of mutually orthogonal minimal tripotents in $C$ (cf. \cite[page 200]{Kaup_ManMath_1997}).\smallskip

Let $F$ be a JBW$^*$-subtriple of a JBW$^*$-triple $E$, and let $e$ be a tripotent in $F$. Clearly, $e$ being minimal in $E$ implies that $e$ is minimal in $F$. We cannot, in general, conclude that $e$ being minimal in $F$ is equivalent to $e$ being minimal in $E$. However, if $F$ is an inner ideal of $E$, Lemma~3.3 in \cite{CabMarPe2024} assures that every minimal tripotent in $F$ is minimal in $E$.\smallskip

For later purposes we need to revisit some structure results and characterizations of atomic JBW$^*$-triples, which are not explicit in the literature. 

\begin{prop}\label{p characterization atomic} Let $E$ be a JBW$^*$-triple. Then the following statements hold:\begin{enumerate}[$(a)$]
		\item $E$ is atomic if, and only if, for every non-zero tripotent $v$ in $E$ there exists a minimal tripotent $e\in E$ satisfying $e\leq v$.
		\item If $E$ is atomic and $e$ is a tripotent in $E$, the Peirce subspaces $E_j (e)$ are atomic JBW$^*$-triples for all $j\in \{0,1,2\}.$ 
	\end{enumerate}  
\end{prop}

\begin{proof}
$(a)$ The ``only if'' implication is clear by structure theory, as we have seen above, every non-zero tripotent $v$ in a Cartan factor, and hence in an atomic JBW$^*$-triple $E$, can be written in the form $v=\hbox{weak$^*$-}\sum_{i} e_i$, where $\{e_i\}_{i}$ is a family of mutually orthogonal minimal tripotents in $E$, and hence $e_i \leq v$ for all $i$.\smallskip

For the ``if'' implication, recall that every element $a$ in a JBW$^*$-triple $E$ can be approximated in norm by a finite linear combination of mutually orthogonal non-zero tripotents $e_1,\ldots, e_k$ in $E$ (cf. \cite[Lemma 3.11]{Horn_MathScand_1987}). Let us fix a non-zero tripotent $e_j$. By hypothesis, there exists a minimal tripotent $v$ satisfying $v\leq e_j$. We can find, via Zorn's lemma, a maximal family $\{v_k^{j}\}_{k}$ of mutually orthogonal minimal tripotents in $E$ with $v_k^{j}\leq e_j$. The series $\sum_{k} v_k^{j}$ is summable in the weak$^*$-topology of $E$ and w$^*$-$\sum_{k} v_k^{j}\leq e_j$ \cite[Corollary 3.13]{Horn_MathScand_1987}. If $e_j - \sum_{k} v_k^{j} \neq 0$, then by the assumptions on $E$, there exists another minimal tripotent $w$ with $w\leq e_j - \sum_{k} v_k^{j}$, which implies that $\{v_k^{j}\}_{k}\cup \{w\}$ is another family of mutually orthogonal minimal tripotents in $E$ bounded by $e_j$, which is impossible. Therefore, $e_j$ can be approximated in the weak$^*$-topology by a finite sum of mutually orthogonal minimal tripotents in $E$. Consequently, $a$ can be approximated in the weak$^*$-topology of $E$ by a finite linear combination of mutually orthogonal minimal tripotents in $E$. \smallskip

$(b)$ By Peirce arithmetic, the Peirce subspaces $E_0(e)$ and $E_2(e)$ are inner ideals of $E$. By applying that $E$ is atomic and $(a)$, we deduce that for each non-zero tripotent $v$ in $E_j (e),$ there exists a minimal tripotent $w\in E$ with $w\leq v$. Observe that in this case $w =\{v,w,v\} \in \{E_j(e),E,E_j(e)\}$, and having in mind that $E_j (e)$ is an inner ideal for all $j =0,2$, we conclude that $w$ is a minimal tripotent in $E_j (e)$ ($j =0,2$). It follows from $(a)$ that $E_j (e)$ is an atomic JBW$^*$-triple for all $j=0,2$.\smallskip

We deal next with $E_1(e)$. By assumptions, $E$ is an atomic JBW$^*$-triple, and hence, $\displaystyle E = \bigoplus_{k}^{\ell_{\infty}} C_k$ for a certain family of Cartan factors $\{C_k\}$. Let $e = (e_k)$ be a tripotent in $E$ (where each $e_k$ is a tripotent in $C_k$). Since $E_j (e)$ coincides with the $\ell_{\infty}$-sum of the family $\{(C_k)_{j} (e_k)\}_{k}$, we can clearly assume that $E$ is a Cartan factor. If $E$ has finite rank, the Peirce-$j$ subspace  $E_j (e)$ also has finite rank, and it is known that in such a case $E_j (e)$ must be a reflexive atomic JB$^*$-triple (\cite[Proposition 4.5]{BuChu92}, \cite[Proof of Theorem 2.3]{BeLoPeRo} and \cite{ChuIochum1987}). We can therefore assume that $E$ is a Cartan factor with infinite rank.\smallskip

We shall distinguish the three remaining cases for $j= 1$.\smallskip

$E=B(H,K)$ is a type 1 Cartan factor, where $H$ and $K$ are infinite dimensional complex Hilbert spaces. Each tripotent $e$ in $E$ is a partial isometry (equivalently, $ee^*$ is a projection in $B(K)$ and $e^*e$ is a projection in $B(H)$), and  $E_1 (e) = ee^* B(H,K) (\mathbf{1}-e^*e) \oplus (\mathbf{1}-ee^*) B(H,K) e^*e$. It is easy to check that the summands $ee^* B(H,K) (\mathbf{1}-e^*e) =  B(ee^*(H),(\mathbf{1}-e^*e) (K))$ are $ B((\mathbf{1}-ee^*)(H),e^*e(K)) $ are orthogonal type 1 Cartan factors, which proves that $E_1 (e)$ is an atomic JBW$^*$-triple. Observe that $E_1(e)$ is not, in general, a factor.\smallskip

Let $j$ be a conjugation on a complex Hilbert space $H$, and set $a^{t} = j a^* j$ for all $a\in B(H)$. Pick a tripotent $e$ in the type 2 Cartan factor $$C= B(H)_{skew} = \{a\in B(H): a^t =- a\}.$$ Clearly, $e$ is a partial isometry in $B(H)$ satisfying $e^t =-e$, $e^* j = - j e$, and $j e^* = - e j$.  For each $a\in C_1 (e)$, the elements $a_1= ee^* a (\mathbf{1}-e^*e)$ and $a_2 =(\mathbf{1}-ee^*) a e^*e$ are orthogonal with $a_2 = -a_1^t$ and $a = a_1 + a_2$.  The mapping $a\mapsto a_1 = ee^* a (\mathbf{1}-e^*e)$ is an isometric triple isomorphism from $C_1 (e)$ onto $B((\mathbf{1}-e^*e)(H),ee^*(H))$, which proves that $C_1(e)$ is a Cartan factor. The case of the type 3 Cartan factor $C= B(H)_{symm} = \{a\in B(H): a^t = a\}$ follows via similar arguments.          
\end{proof}

The following identity principle is almost explicit in \cite[Theorem 3.7]{GarLiPeSu}, we include here a simplified version for completeness reasons.

\begin{prop}\label{p identity principle for truncation preservers} Let $\Delta : E\to E$ be a (non-necessarily additive) bijection preserving truncations and inner quadratic annihilators of elements, where $E$ is an atomic JBW$^*$-triple. Suppose additionally that for each minimal tripotent $e$ in $E$ and each complex number $\lambda$ we have $\Delta (\lambda e ) = \lambda e$. Then $\Delta$ is the identity mapping on $E$.  
\end{prop} 

\begin{proof} Let us take an element $a\in E$, and a minimal tripotent $e\in E$ supported at a pure atom $\varphi_e$ in $E_*$. We reduce to the following two cases:\smallskip
	
\textit{Case 1.} If $\varphi_e(a)\neq 0$ then $\varphi_e(a) e$ is a truncation of $a$ since $$\begin{aligned}  \{\varphi_e(a) e, \varphi_e(a) e, \varphi_e(a) e\} &= \varphi_e(a)^2 \overline{\varphi_e(a)} e = \varphi_e(a)^2 Q(e) P_2(e) (a)  \\
&= \varphi_e(a)^2 Q(e) (a) = \{\varphi_e(a) e, a, \varphi_e(a) e \}.
\end{aligned}$$ The hypotheses on $\Delta$ imply that $\Delta (\varphi_e(a) e) = \varphi_e(a) e$ is a truncation of $\Delta(a),$ and hence $$\begin{aligned} \varphi_e(a)^2 \overline{\varphi_e(a)} e & =
		\{\Delta(\varphi_e(a)e),\Delta(\varphi_e(a)e),\Delta(\varphi_e(a)e)\} \\ 
		&=\{\Delta(\varphi_e(a)e), \Delta(a),\Delta(\varphi_e(a)e)\} = \varphi_e(a)^2 \{e, \Delta(a), e\} \\ 
		&= \varphi_e(a)^2 \overline{\varphi_e(\Delta(a))} e,
	\end{aligned} $$  which gives 
	$$\varphi_e (a) =\varphi_e(\Delta(a)),\hbox{ equivalently, }\varphi_e(a-\Delta(a))=0.$$

	\textit{Case 2.} If $\varphi_e (a)= 0,$ equivalently, $P_2(e)(a)=0 \Leftrightarrow a\in \leftindex^{\perp_{q}}\{e\}$, and thus, the assumptions on $\Delta$ assure that $\Delta (a) \in \leftindex^{\perp_{q}}\{\Delta(e)\} = \leftindex^{\perp_{q}}\{e\}$, therefore  $\varphi_e(\Delta(a)-a) = \varphi_e(\Delta(a)) -\varphi_e(a)  =0$.\smallskip
	
We have proved that $\varphi_e(a-\Delta(a))=0,$ for all $e\in \mathcal{U}_{min}(E),$ $a\in E$. The desired conclusion follows from the fact that pure atoms of $E$ separate the points in $E$.	
\end{proof}

In our next result we describe the Peirce-$1$ subspace associated with a minimal tripotent in an atomic JB$^*$-triple in terms of inner quadratic annihilators. 

\begin{lem}\label{l Peirce1 in annihilators} Let $e$ be a minimal tripotent in an atomic JBW$^*$-triple $E$. Then we have 
	$$E_1 (e) =\leftindex^{\perp_{q}}\{e\}\bigcap \left(\bigcap_{\substack{v\in \mathcal{U}_{min}(E)\\ v\perp e}} \leftindex^{\perp_{q}}\{v\}\right) =\leftindex^{\perp_{q}}\{\mathbb{C} e\}\bigcap \left(\bigcap_{\substack{v\in \mathcal{U}_{min}(E)\\ v\perp e}} \leftindex^{\perp_{q}}\{\mathbb{C} v\}\right).$$
\end{lem}	

\begin{proof} The second equality is clear, we shall only prove the first one.\smallskip 
	
$(\supseteq)$ Take an element $a$ in the intersection given in the right-hand-side above. Since $a\in \leftindex^{\perp_{q}}\{e\}$ it follows that $a\in E_1(e)\oplus E_0(e)$. Similarly $a\in   E_1(v)\oplus E_0(v)$ for all $v\in \mathcal{U}_{min}(E)$ with $v\perp e$, equivalently, for all $v\in \mathcal{U}_{min}(E_0(e))$. Therefore $P_2(v) (P_0(e)(a)) =0$ for all $v\in \mathcal{U}_{min}(E_0(e))$. Having in mind that $E_0(e)$ is an atomic JBW$^*$-triple (see \Cref{p characterization atomic}$(b)$), and the fact that the set $\{P_2 (v) : v\in \mathcal{U}_{min}(E_0(e))\}$ separates the points of $E_0(e)$, we get $P_0(e) (a) =0$, and hence $a = P_1(e) (a) \in E_1 (e)$.\smallskip

$(\subseteq)$ Take now $a\in E_1 (e)$. Clearly $a\in \leftindex^{\perp_{q}}\{e\} = E_0(e)\oplus E_1(e)$. If $v\in \mathcal{U}_{min} (E)$ with $v\perp e$, it follows from \Cref{l orthogonal tripotents do not admit a collinear to one of them but not to the other} that $E_1 (e)\subseteq \leftindex^{\perp_{q}}\{v\}$. 
\end{proof}

We can now show that every surjective additive mapping preserving truncations in both directions between atomic JBW$^*$-triples preserves Peirce-$1$ subspaces associated to minimal tripotents. 

\begin{prop}\label{p preservation of Peirce1} Let $E$ and $F$ be atomic JBW$^*$-triples, and let $A: E\to F$ be a surjective additive mapping preserving truncations in both directions. Then for each minimal tripotent $e$ in $E$ we have 
$ A(E_1 (e)) = F_1 (r(A(e))).$
\end{prop}

\begin{proof} \Cref{l basic properties of additive preservers of truncations into abd images of min trip}$(c)$, $A(\mathbb{C}\ \mathcal{U}_{min} (E)) = \mathbb{C}\ \mathcal{U}_{min} (F)$. Take any $v\in \ \mathcal{U}_{min} (E)$ with $v\perp e$. \Cref{l basic properties of additive preservers of truncations into abd images of min trip}$(d)$ assures that $A(\mathbb{C} v)  = \mathbb{C} A(v) = \mathbb{C} r(A(v))$, while \Cref{prop preservation of orthogonality} proves that $r(A(v))\perp r(A(e))$. This proves that $$ A\left(\mathbb{C}\ \left\{ v\in \mathcal{U}_{min} (E) : v\perp e \right\}\right) \subseteq \mathbb{C}\ \left\{ w\in \mathcal{U}_{min} (F) : w\perp r(A(e)) \right\}.$$ Similarly, since $r(A^{-1} (r(A(e))))\in \mathbb{T} e$ (see \Cref{l basic properties of additive preservers of truncations into abd images of min trip}$(e)$), 
	$$ A^{-1} \left(\mathbb{C}\ \left\{ w\in \mathcal{U}_{min} (F) : w\perp r(A(e)) \right\}\right) \subseteq \mathbb{C}\ \left\{ v\in \mathcal{U}_{min} (E) : v\perp e \right\},$$ and thus \begin{equation}\label{eq image under A of mutliples of minimal trip orthogonal to e}  A\left(\mathbb{C}\ \left\{ v\in \mathcal{U}_{min} (E) : v\perp e \right\}\right) = \mathbb{C}\ \left\{ w\in \mathcal{U}_{min} (F) : w\perp r(A(e)) \right\}. 
	\end{equation} Finally, by \Cref{l Peirce1 in annihilators}, \eqref{eq image under A of mutliples of minimal trip orthogonal to e}, and the fact that $A$ preserves inner quadratic annihilators we derive that 
\begin{align*}
	 A(E_1 (e))  &= A\left( \leftindex^{\perp_{q}}\{\mathbb{C} e\}\bigcap \left(\bigcap_{\substack{v\in \mathcal{U}_{min}(E)\\ v\perp e}} \leftindex^{\perp_{q}}\{\mathbb{C} v\}\right) \right) 	  \end{align*} 
	 \begin{align*}
	  &\subseteq  A\left(\leftindex^{\perp_{q}}\{\mathbb{C} e\}\right) \bigcap \left(\bigcap_{\substack{v\in \mathcal{U}_{min}(E)\\ v\perp e}} A\left(\leftindex^{\perp_{q}}\{\mathbb{C} v\}\right) \right) 	  \end{align*} 
	  \begin{align*} \hspace*{7mm}
	  &=  \leftindex^{\perp_{q}}\{\mathbb{C} r(A(e))\} \bigcap \left(\bigcap_{\substack{v\in \mathcal{U}_{min}(E)\\ v\perp e}} \leftindex^{\perp_{q}}\{\mathbb{C} r(A\left(v\right))\} \right)
	  \end{align*} 
	  \begin{align*} \hspace*{20mm}
	 &= \leftindex^{\perp_{q}}\{\mathbb{C} r(A(e))\}\bigcap \left(\bigcap_{\substack{w\in \mathcal{U}_{min}(F)\\ w\perp r(A(e))}} \leftindex^{\perp_{q}}\{\mathbb{C} w\}\right) = F_1 (r(A(e))). 
	\end{align*} We have therefore shown that $A(E_1 (e)) \subseteq F_1 (r(A(e)))$. Similarly, $$A^{-1}\Big(F_1 \big(r(A(e))\big)\Big) \subseteq E_1 \Big(r\big(A^{-1} (r(A(e))) \big)\Big) = E_1 (e),$$ which gives $A(E_1 (e)) = F_1 (r(A(e)))$.
\end{proof}

\begin{cor}\label{c preservation of Peirce2 for rank2 trip} Let $E$ and $F$ be atomic JBW$^*$-triples, and let $A: E\to F$ be a surjective additive mapping preserving truncations in both directions. Then for each rank-$2$ tripotent $e$ in $E$ we have 
	$ A(E_2 (e)) = F_2 (r(A(e))).$
\end{cor}

\begin{proof} Since $e$ has rank-$2$, we can find two orthogonal minimal tripotents $e_1,e_2$ in $E$ such that $e = e_1 + e_2$. It is well-known that $$E_2 (e) = E_2 (e_1) \oplus E_1(e_1)\cap E_1 (e_2) \oplus E_2(e_2)=\mathbb{C} e_1 \oplus E_1(e_1)\cap E_1 (e_2) \oplus \mathbb{C} e_2.$$ Observe that $A(e) = A(e_1) + A(e_2) = \boldsymbol{\gamma}(e_1,A) r(A(e_1)) + \boldsymbol{\gamma}(e_2,A) r(A(e_2)),$ where $r(A(e_1))$ and $r(A(e_2))$ are two orthogonal minimal tripotents in $F$, and hence $r(A(e)) = r(A(e_1)) + r(A(e_2))$. Furthermore, $A^{-1} \Big( r(A(e)) \Big) = A^{-1} \Big( r(A(e_1)) \Big)+ A^{-1} \Big( r(A(e_2)) \Big)$, and hence $$r\Big(A^{-1} \Big( r(A(e)) \Big)\Big) = \gamma_1 e_1 + \gamma_2 e_2,$$ for suitable $\gamma_1,\gamma_2\in \mathbb{T}$. \smallskip
	
The additivity of $A$, \Cref{l basic properties of additive preservers of truncations into abd images of min trip}, \Cref{prop preservation of orthogonality}, and \Cref{p preservation of Peirce1} lead to 
	$$\begin{aligned}
		  A(E_2 (e)) &= \mathbb{C} r(A(e_1)) \oplus A(E_1(e_1)\cap E_1 (e_2)) \oplus \mathbb{C} r(A(e_2)) \\
		  &\subseteq  \mathbb{C} r(A(e_1)) \oplus A(E_1(e_1))\cap A(E_1 (e_2)) \oplus \mathbb{C} r(A(e_2)) \\
		  &= \mathbb{C} r(A(e_1)) \oplus F_1(r(A(e_1)))\cap F_1 (r(A(e_2))) \oplus \mathbb{C} r(A(e_2)) = F_2 (r(A(e))).
	\end{aligned}$$ That is, $ A(E_2 (e))\subseteq F_2 (r(A(e)))$. Similarly $$A^{-1}\Big( F_2 \big(r(A(e))\big) \Big) \subseteq E_2 \Big(r\big(A^{-1} (r(A(e)) \big)\Big) = E_2 (\gamma_1 e_1 + \gamma_2 e_2) = E_2(e),$$ which concludes the proof. 
\end{proof}

We can now extend some useful consequences of our previous results. 

\begin{prop}\label{c each Cartan factor is mapped to a Cartan factor Badalin-1} Let $E$ and $F$ be atomic JBW$^*$-triples, that is, $\displaystyle E = \bigoplus_{k\in \Gamma_1}^{\ell_{\infty}} C_k$ and $\displaystyle F = \bigoplus_{{j\in \Gamma_2}}^{\ell_{\infty}} \widetilde{C}_j$, where $\{C_k\}_{k\in \Gamma_1}$ and $\{\widetilde{C}_j\}_{j\in \Gamma_2}$ are two families of Cartan factors. Suppose, additionally, that $E$ contains no one-dimensional Cartan factors as direct summands {\rm(}i.e., dim$(C_k)\geq 2$ for all $k${\rm)}. Let $A: E\to F$ be a surjective additive mapping preserving truncations in both directions. Then the following statements hold:\begin{enumerate}[$(a)$]\item $F$ contains no one-dimensional Cartan factors as direct summands.
\item For each $k\in \Gamma_1$ there exists a unique $\sigma (k)\in \Gamma_2$ such that $A(C_k) = \tilde{C}_{\sigma(k)}$, and the mapping $\sigma: \Gamma_1 \to \Gamma_2$ is a bijection. Furthermore, the ranks of $C_k$ and $\tilde{C}_{\sigma(k)}$ coincide.
\end{enumerate}
\end{prop}

\begin{proof} $(a)$ Arguing by contradiction we assume that $\widetilde{C}_{j_0} = \mathbb{C}$ for some $j_0$. Let $e_{j_0}$ be any norm-one element in $\widetilde{C}_{j_0}$. Observe that $e_{j_0}$ is a minimal tripotent in $F$ and $F_1 (e_{j_0}) =\{0\}.$ Let us write $A^{-1} (e_{j_0}) = \boldsymbol{\gamma} (e_{j_0},A^{-1}) r(A^{-1} (e_{j_0})),$ where $r(A^{-1} (e_{j_0}))$ is a minimal tripotent in $E$ and $\boldsymbol{\gamma} (e_{j_0},A^{-1})>0$ (cf. \Cref{l basic properties of additive preservers of truncations into abd images of min trip} and \Cref{r existence of gammaAe and gamm1Ar}). It follows that $r(A^{-1} (e_{j_0}))$ belongs to a unique Cartan factor $C_{k_0}$ in the decomposition of $E$. We know from the hypotheses that dim $(C_{k_0})\geq 2$. Therefore, $E_1 \Big(r\big(A^{-1} (e_{j_0})\big)\Big)\neq \{0\}$. \Cref{l basic properties of additive preservers of truncations into abd images of min trip}$(e)$ implies that $r(A(r(A^{-1} (e_{j_0}))))\in \mathbb{T} e_{j_0}$, and hence \Cref{p preservation of Peirce1} implies that $$ \{0\} \neq A \Big(E_1 \big(r(A^{-1} (e_{j_0}))\big)\Big) = F_1 \Big(r\Big(A\big(r(A^{-1} (e_{j_0}))\big)\Big)\Big) = F_1 (e_{j_0}) = \{0\},$$ which is impossible.\smallskip
		
$(b)$ Fix $k\in \Gamma_1$ and a minimal tripotent $e_k\in C_k$. \Cref{l basic properties of additive preservers of truncations into abd images of min trip}$(c)$ (see also \Cref{r existence of gammaAe and gamm1Ar}) assures that $r(A(e_k))$ is a minimal tripotent in $F$, and hence it belongs to a unique Cartan factor in the decomposition of $F$, we denote this latter factor by $\widetilde{C}_{\sigma(k)}$.\smallskip

Given a minimal tripotent $w\in C_k$ with $w\top e_k,$ \Cref{prop preservation of collinearity} assures that $r(A(w))\top r(A(e_k))$, and thus $r(A(w)), A(w)\in  \widetilde{C}_{\sigma(k)}.$\smallskip

Let $v$ be any other minimal tripotent in $C_k$.\smallskip 

If $C_k$ has rank-one, it is a complex Hilbert space with dimension $\geq 2$, it follows form the discussion in \Cref{r quadratic annihilator and Euclidean orthogonality} that we can find a minimal tripotent $w\in C_k$ such that $w\top e_k$ and $v$ is a linear combination of $w$ and $e_k$. It follows that $A(v) \in A(\mathbb{C} w) + A(\mathbb{C} e_k) =\mathbb{C} r(A( w)) + \mathbb{C} r(A(e_k)) \in \widetilde{C}_{\sigma(k)}$ (cf. \Cref{l basic properties of additive preservers of truncations into abd images of min trip}$(e)$).\smallskip

Assuming that $C_k$ has rank $\geq 2$, it follows from \cite[Lemma 3.10]{Polo_Peralta_AdvMath_2018} that one of the next statements holds:\smallskip

$(1)$ There exists a quadrangle of minimal tripotents $(e_{k},e_2,e_3,e_4)$ such that $v$ is a linear combination of $e_{k},e_2,e_3,$ and $e_4$, and thus $$A(v)\in\mathbb{C} r(A(e_{k})) + \mathbb{C} r(A(e_2)) + \mathbb{C} r(A(e_3))+ \mathbb{C} r(A(e_4))\in  \widetilde{C}_{\sigma(k)}.$$ 

$(2)$ There exists a trangle of the form $(e_k,u,\tilde{e}_k),$ such that $\tilde{e}_k$ is minimal, $u$ is a rank-2 tripotent, and $v$ is a linear combination of $e_k,u,$ and $\tilde{e}_k$. Observe that $\frac{e_{k}\pm u+\tilde{e}_k}{2}$ are minimal tripotents with $\frac{e_{k}\pm u+\tilde{e}_k}{2}\notin \leftindex^{\perp_{q}}\{\mathbb{C} e_k\}$ (cf. \Cref{r linear combinatons of trangles and quadrangles}). It follows that $A(\frac{e_{k}\pm u+\tilde{e}_k}{2}) \notin \leftindex^{\perp_{q}}\{\mathbb{C} A(e_k)\} = \leftindex^{\perp_{q}}\{r(A(e_k))\},$ and thus $A(\frac{e_{k}\pm u+\tilde{e}_k}{2})\in \widetilde{C}_{\sigma(k)}$. We therefore have $$\begin{aligned}
A(v)&\in A(\mathbb{C} e_{k}) + A(\mathbb{C} u) + A(\mathbb{C} \tilde{e}_{k})\\ &\subseteq \mathbb{C}  A(e_{k}) + \mathbb{C}  A\left(\frac{e_{k}+ u+\tilde{e}_k}{2} \right) + \mathbb{C}  A\left(\frac{e_{k}- u+\tilde{e}_k}{2} \right) + \mathbb{C} A( \tilde{e}_{k})\in \widetilde{C}_{\sigma(k)}.
\end{aligned}$$ \smallskip 

The previous arguments show that $$ A\left( \mathbb{C}\  \mathcal{U}_{min} (C_k) \right) \subseteq \mathbb{C} \  \mathcal{U}_{min}(\widetilde{C}_{\sigma(k)}) \subset \widetilde{C}_{\sigma(k)}.$$ By applying a similar reasoning to $A^{-1}$ we deduce that \begin{equation}\label{eq A maps multiples of the mtrip in a factor to multiples of mtrip in a factor} A\left( \mathbb{C}\  \mathcal{U}_{min} (C_k) \right) = \mathbb{C} \  \mathcal{U}_{min}(\widetilde{C}_{\sigma(k)}). 
\end{equation}	 The $\sigma: \Gamma_1 \to \Gamma_2,$  $k\mapsto \sigma(k)$ is well-defined and a bijection because $A$ is bijective.\smallskip

Fix $k_0\in \Gamma_1$, and set $\displaystyle C_{k_0}^{\perp} := \bigoplus_{k\in \Gamma_1, k\neq k_0}^{\ell_{\infty}} C_k$ and $\displaystyle \widetilde{C}_{\sigma(k_0)}^{\perp} = \bigoplus_{{j\in \Gamma_2, j\neq \sigma(k_0)}}^{\ell_{\infty}} \widetilde{C}_j$. Clearly, $C_{k_0}^{\perp}$ and $\widetilde{C}_{\sigma(k_0)}^{\perp}$ are atomic JBW$^*$-triples with $\displaystyle\mathcal{U}_{min} \left( C_{k_0}^{\perp} \right) = \bigcup_{k\in \Gamma_1, k\neq k_0} \mathcal{U}_{min}(C_k)$ and $\displaystyle\mathcal{U}_{min} \left( \tilde{C}_{k_0}^{\perp} \right) = \bigcup_{j\in \Gamma_2, j\neq \sigma(k_0)} \mathcal{U}_{min}(\tilde{C}_j)$. It follows from \eqref{eq A maps multiples of the mtrip in a factor to multiples of mtrip in a factor} that $$ A\left( \mathbb{C}\ \displaystyle\mathcal{U}_{min} \left( C_{k_0}^{\perp} \right) \right) = \mathbb{C}\ \mathcal{U}_{min} \left( \tilde{C}_{\sigma(k_0)}^{\perp} \right).$$ For each $x\in C_{k_0}$, we have $x\perp \mathbb{C} \ \mathcal{U}_{min} \left( C_{k_0}^{\perp} \right)$, which assures that $x$ lies in $\leftindex^{\perp_{q}}{\Big\{\mathbb{C} \ \mathcal{U}_{min} \left( C_{k_0}^{\perp} \right)\Big\}},$ and hence $$A(x)\in \leftindex^{\perp_{q}}{\Big\{A\left(\mathbb{C} \ \mathcal{U}_{min} \left( C_{k_0}^{\perp} \right)\right)\Big\}} = \leftindex^{\perp_{q}}{\Big\{\mathbb{C}\ \mathcal{U}_{min} \left( \tilde{C}_{\sigma(k_0)}^{\perp} \right) \Big\}}.$$ This implies that $P_2 (v) (A(x)) = 0,$ for all $v\in \mathcal{U}_{min} \left( \tilde{C}_{\sigma(k_0)}^{\perp} \right)$. Now, by applying that the set $\left\{P_2(v) : v\in \mathcal{U}_{min} \left( \tilde{C}_{\sigma(k_0)}^{\perp} \right) \right\}$ separates the points in $\tilde{C}_{\sigma(k_0)}^{\perp}$, we arrive to $A(x)\in \tilde{C}_{\sigma(k_0)}$. We have therefore shown that $A\left(C_{k_0}\right) \subseteq \tilde{C}_{\sigma(k_0)},$ and consequently $A\left(C_{k_0}\right) = \tilde{C}_{\sigma(k_0)}.$\smallskip

The final statement is a consequence of \Cref{prop preservation of orthogonality} and the fact that in an atomic JBW$^*$-triple the rank coincides with the cardinality of a maximal set of mutually orthogonal minimal tripotents. 
\end{proof}
 
The next proposition relies on the main result in \cite{LiLiuPe24} combined with the presence of JB$^*$-subtriples isometrically isomorphic to Hilbert spaces of dimension bigger than or equal to $2$. 

\begin{prop}\label{p subtriple generated by two collinear minimal tripotents Badalin-2} Let $E$ and $F$ be JBW$^*$-triples, and let $A: E\to F$ be a surjective additive mapping preserving truncations in both directions. Suppose that $e_1,\ldots, e_n$ are mutually collinear minimal tripotents in $E$ with $n\geq 2$. Then the restriction of $A$ to the JB$^*$-subtriple of $E$ generated by $e_1,\ldots, e_n$ is a positive scalar multiple of a {\rm\!(}complex{)} linear or a conjugate-linear isometry. 
\end{prop}

\begin{proof} By \Cref{l basic properties of additive preservers of truncations into abd images of min trip} and \Cref{prop preservation of collinearity}, the elements $r(A(e_1)),\ldots, r(A(e_n))$ are mutually collinear minimal tripotents in $F$. Observe that the JB$^*$-subtriple generated by $e_1,\ldots, e_n$  (respectively, by $r(A(e_1)),\ldots, r(A(e_n))$) is isometrically isomorphic to an $n$-dimensional complex Hilbert space in which $e_1,\ldots, e_n$ (respectively, $r(A(e_1)),\ldots, r(A(e_n))$) is an orthonormal basis (cf. \cite[Lemma in page 306]{Dang_Friedm_MathScand_1987}). Let $E_0$ and $F_0$ denote the JB$^*$-subtriples of $E$ and $F$ generated by $e_1,\ldots, e_n$ and $r(A(e_1)),\ldots, r(A(e_n)),$  respectively.  A new application of \Cref{l basic properties of additive preservers of truncations into abd images of min trip} and the additivity of $A$ implies that $$ A(E_0) = A(\mathbb{C} e_1 \oplus \ldots \oplus \mathbb{C} e_n ) = \mathbb{C} r(A(e_1)) \oplus \ldots \oplus \mathbb{C} r(A(e_n)) = F_0,$$ that is, $A|_{E_0}$ is a bijective additive mapping preserving truncations in both directions. \Cref{p trunc preservers between Hilbert spaces} gives the desired conclusion.
\end{proof}

Back to the setting of atomic JBW$^*$-triples set in \Cref{c each Cartan factor is mapped to a Cartan factor Badalin-1} we can improve our previous conclusions.   
  
\begin{prop}\label{c each Cartan factor is mapped to a Cartan factor 2 Badalin-3} Under the hypotheses of \Cref{c each Cartan factor is mapped to a Cartan factor Badalin-1}, the following statements hold:\begin{enumerate}[$(a)$] \item For each minimal tripotent $e\in E$ one of the next statements holds:
		\begin{enumerate}\item[$(a.1)$] $A(\lambda e) = \boldsymbol{\gamma}(e,A) \lambda r(A(e))= \lambda A(e),$ for all $\lambda \in \mathbb{C}$. 
			\item[$(a.2)$] $A(\lambda e) = \boldsymbol{\gamma}(e,A) \overline{\lambda} r(A(e)) = \overline{\lambda} A(e),$ for all $\lambda \in \mathbb{C}$.
		\end{enumerate}			
\item If $e,v\in \mathcal{U}_{min} (E)$ with $e\top v$, we have $\boldsymbol{\gamma}(e,A)= \boldsymbol{\gamma}(v,A).$
\item The mapping $A_r: \mathcal{U}_{min} (E) \to \mathcal{U}_{min} (F)$, $e\mapsto A_r(e) :=r(A(e))$ is a bijective mapping preserving collinearity and orthogonality in both directions.
	\end{enumerate}
\end{prop}

\begin{proof}$(a)$ Let us fix an arbitrary minimal tripotent $e\in E$. Suppose $v$ is another minimal tripotent in $E$ such that $v\top e$ (we can always assume the existence of such minimal tripotent $v$ since dim$(C_k)\geq 2$ for all $k\in \Gamma_1$, compare, for example \Cref{sec: rank-one JB*-triples} and \cite[Lemma 3.10]{Polo_Peralta_AdvMath_2018}). Let $E_0$ denote the JB$^*$-subtriple of $E$ generated by $e$ and $v$. By \Cref{p subtriple generated by two collinear minimal tripotents Badalin-2} the restriction of $A$ to $E_0$ is a positive scalar multiple of a (complex) linear or a conjugate-linear isometry, that is, there exists $\gamma_{_{E_0}}>0$ such that $\gamma_{_{E_0}}^{-1} A|_{_{E_0}}$ is a (complex) linear or a conjugate-linear isometry. Consequently, by \Cref{r existence of gammaAe and gamm1Ar}, one of the next statements holds:  
	\begin{enumerate}\item[$(a.1)$] $A(\lambda e) = \boldsymbol{\gamma}(e,A) \lambda r(A(e)) = \lambda A(e)$ for all $\lambda \in \mathbb{C}$, and the same identity holds when $e$ is replaced by any norm-one element in $E_0$. 
		\item[$(a.2)$] $A(\lambda e) = \boldsymbol{\gamma}(e,A) \overline{\lambda} r(A(e)) = \overline{\lambda} A(e)$ for all $\lambda \in \mathbb{C}$, and the same identity holds when $e$ is replaced by any norm-one element in $E_0$.
	\end{enumerate} 
	
It also follows from the above that $$ \boldsymbol{\gamma}(e,A) =\| A(e)\| = \gamma_{_{E_0}} \|e\| = \gamma_{_{E_0}} \|v\| = \|A(v)\| = \boldsymbol{\gamma}(v,A),$$ which proves the statement in $(b)$.\smallskip   

$(c)$ It follows from \Cref{l basic properties of additive preservers of truncations into abd images of min trip}$(a)$ and $(b)$ that the mapping $A_r$ is well-defined and injective. The surjectivity of $A_r$ is a straightforward consequence of $(a)$ and the commented \Cref{l basic properties of additive preservers of truncations into abd images of min trip}. The remaining properties of $A_r$ follow from \Cref{prop preservation of collinearity} and \Cref{prop preservation of orthogonality}.
\end{proof}

As we shall see later, three dimensional spin factors deserve to be treated independently. It is well-known that each three dimensional Hilbert space can be isometrically identified with the complex Banach space $S_2$ of all $2\times 2$ symmetric complex matrices of the form $\left(\begin{matrix} \alpha & \beta \\
	\beta & \delta 
\end{matrix}\right)$ with $\alpha$, $\beta$ and $\delta$ in $\mathbb{C}$ (see, for example, \cite[proof of Lemma 3.4]{Kal_Peralta_Ann_Math_Phys_2021}).

\begin{prop}\label{p S2 to a Cartan factor Badalin-4} Let $F$ be an atomic JBW$^*$-triple expressed in the form $\displaystyle F = \bigoplus_{{j\in \Gamma_2}}^{\ell_{\infty}} \widetilde{C}_j$, where $\{\widetilde{C}_j\}_{j\in \Gamma_2}$ is a family of Cartan factors. Let $A: S_2\to F$ be a surjective additive mapping preserving truncations in both directions. Then the following statements hold:\begin{enumerate}[$(a)$]\item $F = S_2$.
		\item $A$ is a positive scalar multiple of a surjective (complex) linear or a conjugate-linear isometry. 
	\end{enumerate}
\end{prop}

\begin{proof} $(a)$ It follows from \Cref{c each Cartan factor is mapped to a Cartan factor Badalin-1}$(b)$ that $F$ must be a Cartan factor of rank-two. Let $e_1=\left(\begin{matrix} 1 & 0 \\
		0 & 0 
	\end{matrix}\right)$, $e_2=\left(\begin{matrix} 0 & 0 \\
	0 & 1 
	\end{matrix}\right)$, and $u=\left(\begin{matrix} 0 & 1 \\
	1 & 0 
	\end{matrix}\right)$. By combining that $A$ is surjective with \Cref{c each Cartan factor is mapped to a Cartan factor 2 Badalin-3} we get $$ F = A(\mathbb{C} e_1\oplus \mathbb{C} u \oplus \mathbb{C} e_2) = \mathbb{C} r(A(e_1))\oplus \mathbb{C} A(u) \oplus \mathbb{C} r(A(e_2)),$$ where $r(A(e_1))$ and $r(A(e_2))$ are mutually orthogonal minimal tripotents in $F$, and hence $2\leq$dim$(F)\leq 3$. Comparing \cite[Table in page 210]{Kaup_ManMath_1997} or \cite[Table in page 475]{Kaup83MathAnn}, the unique Cartan factor satisfying these properties is $S_2$, that is, the $3$-dimensional spin factor.\smallskip
	
$(b)$ We keep the notation above. Note that $(S_2)_1 (e_1) = (S_2)_1 (e_2) = \mathbb{C} u$, and thus \Cref{p preservation of Peirce1}, implies that $\mathbb{C} A(u)= A(\mathbb{C} u ) = F_1 (r(A(e_1))) =  F_1 (r(A(e_2)))$. It is known that we can find a linear triple automorphism $\Phi: S_2\to S_2$ satisfying $\Phi (r(A(e_1)) ) = e_1$, $\Phi (r(A(e_2)) ) = e_2$ and $\Phi (r(A(u))) = \mu u$ for some $\mu\in \mathbb{T}$ (see, for example, \cite[Proposition 5.8]{Kaup_ManMath_1997}). The mapping $\Phi A: S_2\to S_2$ satisfies the same hypotheses of $A$. We also know that $\mathbb{C} r(A(u)) = F_1(r(A(u))) \ni A(u) = \alpha r(A(u))$ for some positive $\alpha$. Observe that the elements $v_{\pm}:=\frac{e_1\pm u+e_2 }{2}$ are mutually orthogonal minimal tripotents (cf. \Cref{r linear combinatons of trangles and quadrangles}), and thus by \Cref{l basic properties of additive preservers of truncations into abd images of min trip} and \Cref{prop preservation of orthogonality}, their images under $\Phi A$ must be positive scalar multiples of two mutually orthogonal minimal tripotents, that is, 
$$ \Phi A \left( \frac{e_1\pm u+e_2 }{2} \right)  = \frac{\boldsymbol{\gamma}(e_1,A) e_1\pm \alpha \mu u+ \boldsymbol{\gamma}(e_2,A) e_2 }{2} = \left(\begin{matrix} \boldsymbol{\gamma}(e_1,A) & \pm \alpha \mu  \\
	\pm \alpha \mu & \boldsymbol{\gamma}(e_2,A) 
\end{matrix}\right)$$ are positive scalar multiples of mutually orthogonal minimal tripotents, which assures that  $$ \boldsymbol{\gamma}(e_1,A)  \boldsymbol{\gamma}(e_2,A) = \alpha^2 \mu^2  \Rightarrow \boldsymbol{\gamma}(e_1,A)  \boldsymbol{\gamma}(e_2,A) = \alpha^2 \hbox{ and } \mu^2 =1,$$ $$ \boldsymbol{\gamma}(e_1,A)^2 = \alpha^2 |\mu|^2, \hbox{ and } \boldsymbol{\gamma}(e_2,A)^2 = \alpha^2 |\mu|^2.$$ We have therefore shown that $\boldsymbol{\gamma}(e_1,A)= \alpha= \boldsymbol{\gamma}(e_2,A),$ $\mu = \pm 1$. \smallskip

Having in mind \Cref{c each Cartan factor is mapped to a Cartan factor 2 Badalin-3}$(a)$ one of the next cases holds:\smallskip 

\emph{Case $(I)$} $A(\lambda e_1) =  \lambda A(e_1),$ for all  $\lambda \in \mathbb{C}.$ Again, by the dichotomy given by \Cref{c each Cartan factor is mapped to a Cartan factor 2 Badalin-3}$(a)$ one of the next statements is satisfied:\smallskip

\emph{Case $(I.a)$} $A(\lambda v_{+}) = \lambda A(v_+)$ for all $\lambda \in \mathbb{C}$. The identity $$\lambda \frac{ A(e_1)+ A(u)+ A (e_2) }{2}  =  \lambda A(v_+) = A(\lambda v_{+})  = \frac{ A(\lambda e_1)+ A(\lambda u)+ A (\lambda e_2) }{2},$$ combined with the fact that $A(e_1)\perp A(e_2)$ and $A(u)\in  F_1 (r(A(e_1))) =  F_1 (r(A(e_2))),$ implies that $A(\lambda u) = \lambda A(u)$ and $A(\lambda e_2) =\lambda A (e_2)$ for all $\lambda \in \mathbb{C}$.  Therefore, since every $x\in S_2$ writes uniquely in the form $x = \lambda_1 e_1 + \lambda_2 u + \lambda_3 e_2$ for $\lambda_j\in \mathbb{C}$, we arrive to $$A (x) = \lambda_1 A(e_1) + \lambda_2 A(u) + \lambda_3 A(e_2),$$ which clearly assures that $A$ is (complex) linear. We further know that $$\Phi A (x) = \Phi A \left( \left(\begin{matrix} \lambda_1 & \lambda_2  \\
	\lambda_2 & \lambda_3 
\end{matrix}\right) \right) = \boldsymbol{\gamma}(e_1,A) \left(\begin{matrix} \lambda_1 & \mu \lambda_2  \\
\mu \lambda_2 & \lambda_3 
\end{matrix}\right) \ \ (\lambda_1,\lambda_2,\lambda_3\in \mathbb{C}),$$ and consequently $A$ is a positive scalar multiple of a surjective isometry on $S_2$.\smallskip
 
\emph{Case $(I.b)$} $A(\lambda v_{+}) = \overline{\lambda} A(v_+)$ for all $\lambda \in \mathbb{C}$. In this case, the identity $$\overline{\lambda} \frac{ A(e_1)+ A(u)+ A (e_2) }{2}  =  \overline{\lambda} A(v_+) = A(\lambda v_{+})  = \frac{ A(\lambda e_1)+ A(\lambda u)+ A (\lambda e_2) }{2},$$ is incompatible with the assumptions we made at the beginning of this \emph{Case $(I)$} since $A(e_1)\perp A(e_2)$ and $A(u)\in  F_1 (r(A(e_1))) =  F_1 (r(A(e_2))).$\smallskip

\emph{Case $(II)$} $A(\lambda e_1) =  \overline{\lambda} A(e_1)$ for all $\lambda \in \mathbb{C}$, we deduce, via similar arguments to those given in \emph{Case $(I)$} above that $A$ is positive scalar multiple of a conjugate-linear isometry.	
\end{proof}

In the next lemma we shall analyse when two minimal tripotents $e$ and $v$ in a Cartan factor $C$ with rank $\geq 2$ can be connected via a linear combination of the elements in a trangle but not in a quadrangle. Recall that the weak$^*$-closed ideal generated by any minimal tripotent in a Cartan factor $C$ is the whole $C$. 

\begin{lem}\label{l existence of expressions in trangle but not in quadrangles} Let $e$ and $v$ be two minimal tripotents in a Cartan factor $C$. Assume that $C$ has rank$\geq 2,$ and we cannot find a quadrangle of minimal tripotents $(e,e_2,e_3,e_4)$ in $C$ such that $v$ is a linear combination of the elements in this quadrangle. Then $C$ is a Cartan factor of type $3$, and there exists a trangle of the form $(e,u,\tilde{e})$ with $\tilde{e}$ minimal and $u$ having rank-$2$ such that $v$ is a linear combination of $e,u,$ and $\tilde{e}$, and the JBW$^*$-subtriple of $C$ generated by this trangle coincides with $C_2(e+\tilde{e}) = C_2 (u)$ and is isometrically isomorphic to $S_2$.  
\end{lem}

\begin{proof} Lemma 3.10 in \cite{Polo_Peralta_AdvMath_2018} combined with the hypotheses assure the existence of a trangle of the form $(e,u,\tilde{e})$ with $\tilde{e}$ minimal and $u$ having rank-$2$, such that $v$ is a linear combination of $e,u,$ and $\tilde{e}$.\smallskip	  
	
The tripotent $u\in C_1(e)$ governs $e$, and hence $e \notop u$. By combining this observation with the assumptions on $C$, the Classification Scheme in \cite[page 305]{Dang_Friedm_MathScand_1987} (see also Proposition in page 308 in the same reference) assures that $C$ is a Cartan factor of type 3. The rest is clear from \cite[Lemma 3.10]{Polo_Peralta_AdvMath_2018} and the general form of minimal and rank-$2$ tripotents in a type 3 Cartan factor.  
\end{proof}

The behaviour of a surjective additive mapping preserving truncations in both directions on quadrangles of minimal tripotents is even easier to determine. 

\begin{prop}\label{p Badalin-4 quadrangles} Let $E$ and $F$ be atomic JBW$^*$-triples, that is, $\displaystyle E = \bigoplus_{k\in \Gamma_1}^{\ell_{\infty}} C_k$ and $\displaystyle F = \bigoplus_{{j\in \Gamma_2}}^{\ell_{\infty}} \widetilde{C}_j$, where $\{C_k\}_{k\in \Gamma_1}$ and $\{\widetilde{C}_j\}_{j\in \Gamma_2}$ are two families of Cartan factors. Suppose, additionally, that $E$ contains no one-dimensional Cartan factors as direct summands {\rm(}i.e., dim$(C_k)\geq 2$ for all $k${\rm)}. Let $A: E\to F$ be a surjective additive mapping preserving truncations in both directions. Suppose $(e_1,e_2,e_3,e_4)$ is a quadrangle of minimal tripotents in $E$. Then the restriction of $A$ to the JBW$^*$-subtriple $E_0$ of $E$ generated by the elements in the quadrangle is a positive scalar multiple of a {\rm(}complex{\rm)} linear or a conjugate-linear isometry. 
\end{prop}

\begin{proof} Observe first that since $e_1\top e_2\top e_3 \top e_4$, \Cref{c each Cartan factor is mapped to a Cartan factor 2 Badalin-3}$(b)$ implies that $$\boldsymbol{\gamma}(e_1,A) = \boldsymbol{\gamma}(e_2,A) = \boldsymbol{\gamma}(e_3,A) = \boldsymbol{\gamma}(e_4,A).$$ Observe that $w= \frac{e_1+e_2+e_3+e_4}{2}$ is a minimal tripotent in $E$ (see \Cref{r linear combinatons of trangles and quadrangles}), and thus by \Cref{c each Cartan factor is mapped to a Cartan factor 2 Badalin-3}$(a)$ one of the next cases holds:\smallskip
	
\emph{Case $(I)$} $A(\lambda w) = \lambda A(w)$ for all $\lambda\in \mathbb{C}$. The identity $$\begin{aligned}
\lambda& \frac{ A(e_1)+ A (e_2) + A (e_3)+ A (e_4) }{2}  =  \lambda A(w)\\ &= A(\lambda w)  = \frac{ A(\lambda e_1)+ A (\lambda e_2) + A (\lambda e_3)+ A (\lambda e_4) }{2},
\end{aligned}$$ combined with the fact that $r(A(e_1))\top r(A(e_2))\top r(A(e_3))\top r(A(e_4))$, $r(A(e_1))\perp r(A(e_3))$, and $r(A(e_2))\perp r(A(e_4))$ proves that $ A(\lambda e_j) = \lambda  A(e_j)$ for all $j \in\{1,2,3,4\}$, $\lambda \in \mathbb{C}$. Therefore $$ A\left( \lambda_1 e_1+ \lambda_2 e_2+  \lambda_3 e_3+  \lambda_4 e_4 \right) = \lambda_1  A\left( e_1\right) + \lambda_2  A\left( e_2\right) +  \lambda_3  A\left( e_3\right) +  \lambda_4  A\left(e_4 \right) \ (\lambda_j\in \mathbb{C}),$$ assures that the restriction of $A$ to the linear span of $e_1,e_2,e_3,e_4$ is (complex) linear, while, $$\begin{aligned}
\lambda_1  &A\left( e_1\right) + \lambda_2  A\left( e_2\right) +  \lambda_3  A\left( e_3\right) +  \lambda_4  A\left(e_4 \right) \\ &= \boldsymbol{\gamma}(e_1,A) \left( \lambda_1 r(A\left( e_1\right)) + \lambda_2  r(A\left( e_2\right)) +  \lambda_3  r(A\left( e_3\right)) +  \lambda_4  r(A\left(e_4 \right)) \right)  \ (\lambda_j\in \mathbb{C}),
\end{aligned} $$ implies that the restriction of $A$ to $E_0=\hbox{span}\{e_1,e_2,e_3,e_4\}$ is a positive scalar multiple of an isometry. It can be easily checked that the minimal tripotents $r(A(e_1)),$ $r(A(e_2)),$ $r(A(e_3)),$ $r(A(e_4))$ form a quadrangle in $C_{\sigma(k)} \subseteq F$, and $A(E_0)$ $=$ $\hbox{span} \{r(A(e_1)), r(A(e_2)), r(A(e_3)), r(A(e_4))\}$.\smallskip	 
	
\emph{Case $(II)$} $A(\lambda w) = \overline{\lambda} A(w)$ for all $\lambda\in \mathbb{C}$. Similar arguments to those given above prove that the restriction of $A$ to $E_0=\hbox{span}\{e_1,e_2,e_3,e_4\}$ is a positive scalar multiple of a conjugate-linear isometry, and $A|_{E_0} : E_0\to A(E_0)$ is a positive multiple of a conjugate-linear triple isomorphism.
\end{proof}

\subsection{Connections with preservers of triple transition pseudo-probabilities} \ \smallskip

As we have seen in \Cref{c each Cartan factor is mapped to a Cartan factor 2 Badalin-3} Under the hypotheses of \Cref{c each Cartan factor is mapped to a Cartan factor Badalin-1}, given two atomic JBW$^*$-triples $E$ and $F$, such that $E$ contains no one-dimensional Cartan factor summands, and a surjective additive mapping $A:E\to F$ preserving truncations in both directions, we can associate $A$ with a bijection $A_r: \mathcal{U}_{min} (E) \to \mathcal{U}_{min} (F)$ preserving collinearity and orthogonality in both directions which is given by $A_r(e) :=r(A(e)).$ However, we lack of a result saying that $A$ maps tripotents to positive scalar multiples of tripotents. So, we must restrict ourself to the sets of minimal tripotents. The remaining part of our arguments will rely on the connections between surjective additive mappings preserving truncations and preservers of triple transition pseudo-probabilities. Let $e$ and $v$ be two minimal tripotents in a JBW$^*$-triple $E$. We recall that the \emph{triple transition pseudo-probability} (of transitioning) from $e$ to $v$ is the complex number given by \begin{equation}\label{eq triple transition pseudo prbability} TTP(e,v)= \varphi_v(e),
\end{equation} where $\varphi_v$ is the unique pure atom in $E_*$ supported at $v$ (cf. \cite{Per_RIMS_2023,Peralta_ResMath_2023}). The triple transition pseudo-probability between minimal tripotents is a symmetric mapping, that is, $TTP(e,v)= \overline{TTP(v,e)},$ for all $e,v\in \mathcal{U}_{min} (E).$ In the case that $e$ and $v$ are minimal projections in $B(H)$, the triple transition pseudo-probability from $e$ to $v$ is precisely the usual transition probability from $e$ to $v$ in the celebrated Wigner's theorem. It is known that every (linear) triple isomorphism $T$ between JBW$^*$-triples is automatically weak$^*$-continuous (cf. \cite[Corollary 3.22]{Horn_MathScand_1987}), and hence $T$ preserves triple transition pseudo-probabilities between minimal tripotents. Reciprocally, the main result in \cite{Peralta_ResMath_2023} shows that every bijective mapping preserving triple transition pseudo-probabilities between the sets of minimal tripotents of two atomic JBW$^*$-triples is precisely the restriction of a (complex-)linear triple isomorphism between the corresponding JBW$^*$-triples.\smallskip

\begin{rem}\label{r conjugate-linear triple isomorphisms} Let $S : E\to F$ be a conjugate-linear triple isomorphism between two JBW$^*$-triples. It is also known that in this case $S$ must be also weak$^*$-continuous (see, for example, \cite[Proposition 2.3]{MarPe2000}). Let us take two minimal tripotents $e,v\in E$ and the corresponding pure atoms $\varphi_{e},\varphi_v$ supported at $e$ and $v$, respectively. Clearly $S(e)$ and $S(v)$ are minimal tripotents in $F$ whose supporting functional are denoted by $\varphi_{S(e)},\varphi_{S(v)},$ respectively. By definition we have $$\begin{aligned}
\varphi_{S(v)} (S(e)) S(v)&= P_2 (S(v)) (S(e)) = \{S(v),\{S(v),S(e),S(v)\},S(v)\} \\
&=  \{S(v),S(\{v,e,v\}),S(v)\} =  S(P_2(v) (e)) = S(\varphi_{v} (e) v)\\
& = \overline{\varphi_{v}} (e) S(v),		
\end{aligned}$$ which assures that $$ TTP(S(e),S(v))= \overline{TTP(e,v)} = TTP(v,e).$$
\end{rem}

\begin{rem}\label{r TTP does not change when computed with respect to a JBW*-subtriple} 
Let $F$ be a JBW$^*$-subtriple of a JBW$^*$-triple $E$. Let us consider two tripotents $e,v\in F$ which are minimal tripotents in $E$. Then the triple transition pseudo-probability from $e$ to $v$ does not change when computed in $E$ or in $F$. Namely, for each minimal tripotent $e\in E$ belonging also to $F$ we have $E_2 (e)= \mathbb{C} e = F_2 (e)$. If $\phi$ is any functional in $E_*$ such that $\|\phi\| = 1 = \phi (e) = 1$, we have $\phi = \phi P_2(e) = \phi|_{E_2(e)} P_2 (e)$ (cf. \cite[Proposition 1]{FriedRusso_Predual}). So, if $\varphi^{E}_{v} \in \partial_{e} (\mathcal{B}_{E_*})$ and $\varphi^{F}_{v} \in \partial_{e} (\mathcal{B}_{F_*})$ denote the unique pure atoms of $E$ and $F$ supported at $v$, respectively, it follows that $\varphi^{E}_{v}|_{F} = \varphi^{F}_{v}$ --in other words, there is a unique norm-preserving extension of $\varphi^{F}_{v}$ to an element in $E_*$.
\end{rem}

We are now in a position to culminate our technical arguments.      
  
\begin{prop}\label{c each Cartan factor is mapped to a Cartan factor 2 Badalin-5} Under the hypotheses of \Cref{c each Cartan factor is mapped to a Cartan factor Badalin-1}, the following statements hold:\begin{enumerate}[$(a)$] 
		\item For each $k\in \Gamma_1$ there exists a {\rm(}unique{\rm)} positive constant $\gamma_k$ such that $\gamma_k = \boldsymbol{\gamma}(w,A)$ for all minimal tripotent $w\in C_k$, that is, $\boldsymbol{\gamma}(w,A)$ does not change when $w$ runs in the set, $\mathcal{U}_{min} (C_k),$ of all minimal tripotents of each Cartan factor $C_k$.  
		\item For each $k\in \Gamma_1$ one of the next statements holds:
		\begin{enumerate}[$(b.1)$]\item $A(\lambda e) = \gamma_k \lambda r(A(e))= \lambda A(e),$ for all $e\in  \mathcal{U}_{min} (C_k)$ and $\lambda \in \mathbb{C}$. 
			\item $A(\lambda e) = \gamma_k \overline{\lambda} r(A(e)) = \overline{\lambda} A(e),$ for all $e\in  \mathcal{U}_{min} (C_k)$ and $\lambda \in \mathbb{C}$. 
		\end{enumerate}
		\item For each $k\in \Gamma_1$ one of the next statements holds:
		\begin{enumerate}[$(c.1)$]\item If $A(\lambda w) = \lambda A(w),$ for all $w\in  \mathcal{U}_{min} (C_k)$ and $\lambda \in \mathbb{C}$, the mapping $A_r$ preserves triple transition pseudo-probabilities between elements in $\mathcal{U}_{min} (C_k)$, that is, $TTP(A_r(e),A_r(v))= TTP(e,v)$, for all $e,v\in \mathcal{U}_{min} (C_k)$.
			\item If $A(\lambda w) =  \overline{\lambda} A(w),$ for all $w\in  \mathcal{U}_{min} (C_k)$ and $\lambda \in \mathbb{C}$, the mapping $A_r$ reverses triple transition pseudo-probabilities between elements in $\mathcal{U}_{min} (C_k)$, that is, $TTP(A_r(e),A_r(v))= \overline{TTP(e,v)} = TTP(v,e)$, for all $e,v\in \mathcal{U}_{min} (C_k)$.
		\end{enumerate}
		\item The families $(\gamma_k)_{k\in \Gamma_1}$ and $(\gamma_k^{-1})_{k\in \Gamma_1}$ are bounded.	
	\end{enumerate}
\end{prop}

\begin{proof} Pick an arbitrary $k\in \Gamma_1$ and two minimal tripotents $e,v\in C_k$. We shall distinguish three main cases.\smallskip

\emph{Case $(1)$:} $C_k$ has rank-one. By \Cref{c each Cartan factor is mapped to a Cartan factor Badalin-1}$(b)$, the restricted mapping $A|_{C_k}: C_k \to \tilde{C}_{\sigma(k)}$ is a surjective mapping preserving truncations in both directions between two rank-one Cartan factors with dim$(C_k)\geq 2$. \Cref{p trunc preservers between Hilbert spaces} proves that $A|_{C_k}$ is a positive scalar multiple of a (complex) linear or a conjugate-linear isometry, and hence the desired conclusion in $(a)$ and $(b)$ are trivially true.\smallskip

\emph{Case $(2)$:} $C_k$ has rank $\geq 2$ and there exists a quadrangle of minimal tripotents of the form $(e,e_2,e_3,e_4)$ in $C_k$ such that $v$ is a linear combination of the elements in the quadrangle. We deduce from \Cref{p Badalin-4 quadrangles} that the restriction of $A$ to the JBW$^*$-subtriple $E_0\subseteq C_k$ spanned by the set $\{e,e_2,e_3,e_4\}$ is a positive scalar multiple of a(complex) linear or a conjugate-linear isometry. Consequently, the conclusions in $(a)$ and $(b)$ trivially hold for $e$ and $v$, and the arbitrariness of thee two elements concludes the proof.\smallskip

\emph{Case $(3)$:} $C_k$ has rank $\geq 2$, but $v$ cannot be written as a linear combination of the elements in a quadrangle of minimal tripotents of the form $(e,e_2,e_3,e_4)$ in $C_k$. \Cref{l existence of expressions in trangle but not in quadrangles} implies that $C_k$ is a type 3 Cartan factor and there exists a trangle of the form $(e,u,\tilde{e})$ with $\tilde{e}$ minimal and $u$ having rank-$2$ such that $v$ is a linear combination of $e,u,$ and $\tilde{e}$, and the JBW$^*$-subtriple of $C$ generated by this trangle coincides with $C_2(e+\tilde{e}) = C_2 (u)$ and is isometrically isomorphic to $S_2$. \Cref{c preservation of Peirce2 for rank2 trip} shows that $A(C_2(e+\tilde{e})) = F_2(r(A(e+\tilde{e})))$. We also know that $r(A(e+\tilde{e}))$ is a rank-$2$ tripotent in $F$ (cf. \Cref{l basic properties of additive preservers of truncations into abd images of min trip} and \Cref{prop preservation of orthogonality} combined with the additivity of $A$), and $F_2(r(A(e+\tilde{e})))$ is an atomic JBW$^*$-triple (cf. \Cref{p characterization atomic}). We can now conclude from \Cref{p S2 to a Cartan factor Badalin-4} the $F_2(r(A(e+\tilde{e}))) = S_2$ and the restriction of $A$ to $E_0 =C_2(e+\tilde{e})$ is a positive scalar multiple of a surjective (complex) linear or a conjugate-linear isometry. As before, this shows that $(a)$ and $(b)$ hold for $e$ and $v$. \smallskip  

$(c)$ We have shown in the proof of $(a)$ and $(b)$ that for each $k\in \Gamma_1,$ and any two minimal tripotents $e,v\in C_k$, we can find a JBW$^*$-subtriple $E_0\subseteq C_k$ containing $e$ and $v$ such that $\gamma_k^{-1} A|_{E_0}: E_0 \to A(E_0)$ is a linear or a conjugate-linear isometry, equivalently, a linear or a conjugate-linear triple isomorphism between JBW$^*$-triples (see \cite[Proposition 5.5]{Kaup83MathAnn}). It is also clear from $(a)$ and \Cref{l basic properties of additive preservers of truncations into abd images of min trip} that $A_r|_{\mathcal{U}_{min}(E_0)} = \gamma_k^{-1} A|_{\mathcal{U}_{min}(E_0)}.$ It follows that $TTP(A_r(e),A_r(v))= TTP(e,v)$ if $A(\lambda w) = \gamma_k \lambda r(A(w))= \lambda A(w),$ for all $w\in  \mathcal{U}_{min} (C_k)$ and all $\lambda \in \mathbb{C},$ and $TTP(A_r(e),A_r(v))= \overline{TTP(e,v)} = TTP(v,e)$ if $A(\lambda w) = \gamma_k \overline{\lambda} r(A(w)) = \overline{\lambda} A(w),$ for all $w\in  \mathcal{U}_{min} (C_k)$ and $\lambda \in \mathbb{C}$ (cf. \Cref{r conjugate-linear triple isomorphisms}, the preceding comments to this result, and \Cref{r TTP does not change when computed with respect to a JBW*-subtriple}).\smallskip

$(d)$ For each $k\in \Gamma_1$ pick a minimal tripotent $e_k \in C_k$. The element $e= (e_k)_{k\in\Gamma_1}$ is a tripotent in $E$ satisfying $$ \gamma_{k_0}= \|A(e_{k_0})\| = \|A(e_{k_0}) + A(e-e_{k_0}) \| = \|A(e) \|.$$ If in this argument we replace $A$ with $A^{-1}$ we obtain the boundedness of the family $(\gamma_{k}^{-1})_{k\in \Gamma_1}$. 	
\end{proof}

 We can now state our main result on surjective additive mappings between atomic JBW$^*$-triples preserving truncations in both directions.
 
\begin{thrm}\label{t main theorem} Let $E$ and $F$ be atomic JBW$^*$-triples, that is, $\displaystyle E = \bigoplus_{k\in \Gamma_1}^{\ell_{\infty}} C_k$ and $\displaystyle F = \bigoplus_{{j\in \Gamma_2}}^{\ell_{\infty}} \widetilde{C}_j$, where $\{C_k\}_{k\in \Gamma_1}$ and $\{\widetilde{C}_j\}_{j\in \Gamma_2}$ are two families of Cartan factors. Suppose, additionally, that $E$ contains no one-dimensional Cartan factors as direct summands {\rm(}i.e., dim$(C_k)\geq 2$ for all $k${\rm)}. Let $A: E\to F$ be a surjective additive mapping. Then the following statements are equivalent:  
	\begin{enumerate}[$(a)$]
		\item $A$ preserves truncations in both directions. 
		\item There exists a  bijection $\sigma: \Gamma_1\to \Gamma_2$, a bounded family $(\gamma_k)_{k\in \Gamma_1}\subseteq \mathbb{R}^+$, and a family $(\Phi_k)_{k\in \Gamma_1},$ where each $\Phi_k$ is a {\rm(}complex{\rm)} linear or a conjugate-linear {\rm(}isometric{\rm)} triple isomorphism from $C_k$ onto $\widetilde{C}_{\sigma(k)}$ satisfying $\inf_{k} \{\gamma_k \} >0,$ and $$A(x) = \Big( \gamma_{k} \Phi_k \left(\pi_k(x)\right) \Big)_{k\in\Gamma_1},\  \hbox{ for all } x\in E,$$ where $\pi_k$ denotes the canonical projection of $E$ onto $C_k.$  
	\end{enumerate} Moreover, if $A$ preserves truncations in both directions and there exists an element $x_0$ in $E$ satisfying $\|A\pi_k (x_0)\| = \|\pi_k (x_0)\| \neq 0$, for all $k\in \Gamma_1$, then the mapping $A$ is a real-linear {\rm(}isometric{\rm)} triple isomorphism.
\end{thrm}

\begin{proof} We shall only prove that $(a)\Rightarrow (b)$. Let $\sigma: \Gamma_1\to \Gamma_2$ denote the bijection given by \Cref{c each Cartan factor is mapped to a Cartan factor Badalin-1}$(b)$, and let $(\gamma_k)_{k\in \Gamma_1}$ stand for the bounded family of positive numbers whose existence is assured by \Cref{c each Cartan factor is mapped to a Cartan factor 2 Badalin-5}. \smallskip
	
Fix now $k\in \Gamma_1$. By \Cref{c each Cartan factor is mapped to a Cartan factor Badalin-1} and \Cref{l basic properties of additive preservers of truncations into abd images of min trip}, the restricted mapping $A|_{C_k} : C_k\to \widetilde{C}_{\sigma(k)}$ is a bijection preserving truncations in both directions. \Cref{c each Cartan factor is mapped to a Cartan factor 2 Badalin-5} implies that one (and just one) of the next statements holds:
	\begin{enumerate}[$(1)$]
		\item $A(\lambda e) = \gamma_k \lambda r(A(e))= \lambda A(e)= \gamma_k \lambda A_r(e),$ for all $e\in  \mathcal{U}_{min} (C_k)$ and $\lambda \in \mathbb{C}$, and $TTP(A_r(e),A_r(v))= TTP(e,v)$, for all $e,v\in \mathcal{U}_{min} (C_k)$..
		\item $A(\lambda e) = \gamma_k \overline{\lambda} r(A(e)) = \overline{\lambda} A(e)= \gamma_k \overline{\lambda} A_r(e),$ for all $e\in  \mathcal{U}_{min} (C_k)$ and $\lambda \in \mathbb{C}$, and $TTP(A_r(e),A_r(v))= \overline{TTP(e,v)} = TTP(v,e)$, for all $e,v\in \mathcal{U}_{min} (C_k)$.
	\end{enumerate}
	
	In case $(1)$, the mapping $A_r^{-1} : \mathcal{U}_{min} (\widetilde{C}_{\sigma(k)})\to  \mathcal{U}_{min} (C_k)$ is a bijection preserving triple transition pseudo-probabilities (see also \Cref{c each Cartan factor is mapped to a Cartan factor 2 Badalin-3}$(c)$). Corollary 3.3 in \cite{Peralta_ResMath_2023} assures the existence of a linear (isometric) triple isomorphism $\Psi_k : \widetilde{C}_{\sigma(k)}\to  C_k$ whose restriction to  $\mathcal{U}_{min} (\widetilde{C}_{\sigma(k)})$ is $A_r^{-1}$. We also know from \Cref{c each Cartan factor is mapped to a Cartan factor 2 Badalin-5} and the assumptions in this case $(1)$, that $A|_{\mathcal{U}_{min} (C_k)} = \gamma_k A_r $, and the mapping $T_k = \gamma_k^{-1} \Psi_k  A|_{C_k} : C_k\to C_k$ is an additive bijection preserving truncations in both directions, which also satisfies $T_k (\lambda e) = \lambda e$ for all $e\in \mathcal{U}_{min} (C_k)$ and $\lambda \in \mathbb{C}$. The identity principle in \Cref{p identity principle for truncation preservers} now implies that $T_k$ is the identity mapping on $C_k$, and thus $A|_{C_k} = \gamma_k \Psi_k^{-1}$ is a positive scalar multiple of an isometric linear triple isomorphism $\Phi_k = \Psi_k^{-1}$ from $C_k$ onto $\widetilde{C}_{\sigma(k)}$.\smallskip
	
	Suppose now that $A|_{C_k} : C_k\to \widetilde{C}_{\sigma(k)}$ satisfies the properties in case $(2)$. It is known that we can always find a conjugation (i.e., a conjugate-linear triple automorphism of period-$2$) $\Psi_k : \widetilde{C}_{\sigma(k)}\to \widetilde{C}_{\sigma(k)}$ (cf. \cite[Theorem 4.1]{Kaup_ManMath_1997} and \cite{Loos77Irvine} for additional details in the case of exceptional Cartan factors). The mapping $\Psi_k A|_{C_k} : C_k\to \widetilde{C}_{\sigma(k)}$ satisfies the properties in case $(1),$ and hence, by the previous conclusion $\Psi_k A|_{C_k} = \alpha_k \widehat{\Phi}_k,$ where $\alpha_k\in \mathbb{R}^{+}$ and $\widehat{\Phi}_k : C_k\to \widetilde{C}_{\sigma(k)}$ is a linear and isometric triple isomorphism. Necessarily $\alpha_k = \gamma_k$. Finally, $ A|_{C_k} = \gamma_k \Psi_k \widehat{\Phi}_k = \gamma_k {\Phi}_k,$ where ${\Phi}_k =\Psi_k \widehat{\Phi}_k: C_k\to \widetilde{C}_{\sigma(k)}$ is a conjugate-linear triple isomorphism.\smallskip
	
	The final statement is clear since the existence of such element $x_0\in E$ implies that $\gamma_k=1$ for all $k$. 
\end{proof}

Our last corollary, which is interesting by itself, is a straightforward consequence of our main theorem.

\begin{cor}\label{c main theorem for a single Cartan factor} Let $C$ and $\tilde{C}$ be Cartan factors with dim$(C)\geq 2$. Then every surjective additive mapping $A: C\to \tilde{C}$ preserving truncations in both directions is a positive scalar multiple of a linear or a conjugate-linear (isometric) triple isomorphism. Moreover, if we also assume that $\|A(x_0)\| = \|x_0\|$ for some non-zero $x_0$ in $C$, the mapping $A$ is a linear or a conjugate-linear {\rm(}isometric{\rm)} triple isomorphism.
\end{cor}

\textbf{Acknowledgements} L. Li supported by National Natural Science Foundation of China (Grant No. 12171251). A.M. Peralta supported by grant PID2021-122126NB-C31 funded by MICIU/AEI/10.13039/501100011033 and by ERDF/EU, by Junta de Andalucía grant FQM375, IMAG--Mar{\'i}a de Maeztu grant CEX2020-001105-M/AEI/10.13039/501100011033 and (MOST) Ministry of Science and Technology of China grant G2023125007L. \smallskip

Part of this work was completed during a visit of A.M. Peralta to Nankai University and the Chern Institute of Mathematics, which he thanks for the hospitality.\medskip

\subsection*{Statements and Declarations} 

All authors declare that they have no conflicts of interest to disclose.

\subsection*{Data availability}

There is no data associate for this submission.



\end{document}